\documentclass[11pt]{article}
\usepackage{amsmath,amssymb,amsthm,lscape,stmaryrd}
\title{On Lie induction and the exceptional series}
\author{Jan E. Grabowski}
\date{20th September 2004}

\newtheorem{definition}{Definition}[section]
\newtheorem{theorem}{Theorem}[section]
\newtheorem*{theorem*}{Theorem}
\newtheorem{proposition}[theorem]{Proposition}
\newtheorem{lemma}[theorem]{Lemma}

\newtheorem{definition+}[theorem]{Definition}
\newenvironment{jegnote}[1]{\par \noindent \textit{Note{#1}:}}{}

\newcounter{property}
\newcounter{subproperty}
\newenvironment{property}{\refstepcounter{property} \par \vspace{1em} \noindent (\arabic{property})\begin{math}\quad }{\end{math} \vspace{1em} \\}

\newenvironment{equivproperty}{\par \vspace{1em} \noindent (\arabic{property}$'$)\begin{math} \quad }{\end{math} \vspace{1em} \\}

\input cyracc.def

\newcommand{\ad}[2]{{\mathrm{ad}}_{#1}(#2)}
\newcommand{\bigdsum}{\bigoplus}
\newcommand{\chr}{\mathop{\mathrm{char}}}
\newcommand{\complex}{\ensuremath \mathbb{C}}
\newcommand{\corank}{\mathop{\mathrm{corank}}}
\newcommand{\dbos}[3]{\Lie{#1} \rtimesdot \Lie{#2} \ltimesdot \Lie{#3}}
\newcommand{\defeq}{\stackrel{\scriptscriptstyle{\mathrm{def}}}{=}}
\newcommand{\dform}{\mathrm{d}}
\newcommand{\dsum}{\ensuremath{ \oplus}}
\newcommand{\dual}[1]{\ensuremath {#1}^{*}}
\renewcommand{\epsilon}{\varepsilon}
\newcommand{\height}[1]{\mathop{\mathrm{ht}}( #1 )}
\newcommand{\id}{\ensuremath \mbox{\textup{id}}}
\renewcommand{\iff}{\ensuremath \Longleftrightarrow}

\newcommand{\inj}{\hookrightarrow}
\newcommand{\integ}{\ensuremath{\mathbb{Z}}}
\newcommand{\inv}[1]{\ensuremath {#1}^{-1}}
\newcommand{\ip}[2]{\ensuremath \lgen\;\!#1,#2\;\!\rgen}
\newcommand{\iso}{\ensuremath \cong}
\newcommand{\Ker}[1]{\ensuremath \mbox{Ker}\>{#1}}
\newcommand{\laction}{\triangleright}
\newcommand{\lgen}{\ensuremath \mathopen{<}} %see also \rgen

\newcommand{\ltimesdot}{\mathrel{\cdot\joinrel\mkern-14.3mu\rhd\mkern-8.7mu<}}
\newcommand{\Lbracket}[2]{\ensuremath{ [\, #1 , #2\, ]}}
\newcommand{\Lbbracket}[2]{\ensuremath{ \llbracket\, #1 , #2\, \rrbracket}}
\newcommand{\Lie}[1]{\ensuremath{\mathfrak{#1}}}
\newcommand{\modcat}[1]{\ensuremath \!_{#1}\mathcal{M}}
\newcommand{\mult}[2]{\mathrm{mult}_{#1}( #2 )}
\newcommand{\nat}{\ensuremath \mathbb{N}}
\newcommand{\octo}{\ensuremath \mathbb{O}}
\newcommand{\op}[1]{{#1}^{\mbox{\scriptsize \textup{op}}}}
\renewcommand{\phi}{\varphi}
\newcommand{\qea}[1]{\ensuremath{U_{q}(\Lie{#1})}}
\newcommand{\rank}{\mathop{\mathrm{rank}}}
\newcommand{\rgen}{\ensuremath \mathclose{>}}
\newcommand{\rtimesdot}{\mathrel{>\joinrel\mkern-6mu\lhd\mkern-7.2mu\cdot}}

\newcommand{\setspan}[2]{\mathrm{span}_{#1}\{\, #2 \,\}}
\newcommand{\smun}[2]{{#1}_{\scriptscriptstyle #2}}
\newcommand{\Ss}{\scriptstyle}
\newcommand{\supp}[1]{\mathop{\mathrm{Supp}}( #1 )}
\newcommand{\Sym}[2]{\mathop{\mathrm{Sym}}^{#1}( #2 )}
\newcommand{\tensor}{\ensuremath \otimes}
\newcommand{\union}{\mathrel{\cup}}
\newcommand{\Wedge}[2]{\mathop{\textstyle{\bigwedge}^{#1}}#2 }
\newcommand{\zeroset}{\ensuremath{\{ 0 \}}}

\newcommand{\ie}{i.e.\ }
\newcommand{\eg}{e.g.\ }
\newcommand{\href}[2]{#2}
\newcommand{\webtilde}{\small{$\sim\!$} \normalsize}

\begin{document}
\maketitle

\begin{abstract}
Lie bialgebras occur as the principal objects in the infinitesimalisation of the theory of quantum groups---the semi-classical theory.  Their relationship with the quantum theory has made available some new tools that we can apply to classical questions.  In this paper, we study the simple complex Lie algebras using the double-bosonisation construction of Majid.  This construction expresses algebraically the induction process given by adding and removing nodes in Dynkin diagrams, which we call Lie induction.

We first analyse the deletion of nodes, corresponding to the restriction of adjoint representations to subalgebras.  This uses a natural grading associated to each node.  We give explicit calculations of the module and algebra structures in the case of the deletion of a single node from the Dynkin diagram for a simple Lie (bi-)algebra.  

We next consider the inverse process, namely that of adding nodes, and give some necessary conditions for the simplicity of the induced algebra.  Finally, we apply these to the exceptional series of simple Lie algebras, in the context of finding obstructions to the existence of finite-dimensional simple complex algebras of types $E_9$, $F_5$ and $G_3$.  In particular, our methods give a new point of view on why there cannot exist such an algebra of type $E_9$.
\end{abstract}

\section{Introduction}

The study of Lie algebras is long-established and widely utilized but in recent years, two generalisations of the theory have become prominent and almost as ubiquitous.  These are the ``semi-classical'' theory of Lie bialgebras and the ``quantum'' theory of quantized enveloping algebras.  In this scheme, the original theory is called ``classical''.  This paper is concerned with the implications of a tool developed in the quantum theory for the classical, via the semi-classical.

We recall that a complex Lie algebra $\Lie{g}$ is a $\complex$-vector space equipped with a map $\Lbracket{\ }{\ }: \Lie{g} \tensor \Lie{g} \to \Lie{g}$, called the bracket, which is anti-symmetric and satisfies the Jacobi identity: \[ \Lbracket{x}{\Lbracket{y}{z}} + \Lbracket{y}{\Lbracket{z}{x}} + \Lbracket{z}{\Lbracket{x}{y}} =0. \]  A Lie bialgebra is a Lie algebra with a cobracket structure $\delta: \Lie{g} \to \Lie{g} \tensor \Lie{g}$ satisfying axioms precisely dual to those for a Lie algebra, with an appropriate compatibility condition.  The definition, due to Drinfel\cprime d (\cite{Drinfeld1}), is comparable with that of a Hopf algebra $H$, where we have a multiplication map $m:H\tensor H \to H$ and a compatible comultiplication $\Delta:H \to H\tensor H$.  The comultiplication defines an algebra structure on the dual $\dual{H}$, so we think of a Hopf algebra as being self-dual in this sense.  The definition of a Lie bialgebra is the semi-classical version or \emph{infinitesimalisation} of this.

An important class of Lie bialgebras is that of quasitriangular Lie bialgebras.  Here, the cobracket $\delta$ is of a specific form, namely, the coboundary of an element $r\in \Lie{g}\tensor \Lie{g}$ satisfying two conditions, one of which is the well-known classical Yang-Baxter equation.  We will mostly be concerned with the Lie algebra structures in this paper and we will usually take the canonical Drinfel\cprime d--Sklyanin cobracket, which does give a quasitriangular Lie bialgebra.  However, a quasitriangular bialgebra structure is absolutely essential to the construction we use.

We also need a generalisation of the notion of a Lie bialgebra, that of a braided-Lie bialgebra.  Given a quasitriangular Lie bialgebra $\Lie{g}$, this is a $\Lie{g}$-covariant bialgebra in the category of $\Lie{g}$-modules where the cobracket $\underline{\delta}$ has non-zero coboundary.  In fact $\dform \underline{\delta} = \psi$, where $\psi$ is a canonical braiding operator.  This definition has been given by Majid (\cite{BraidedLie}) as the infinitesimal version of his braided groups for Hopf algebras.  Full details of these definitions may be found in Section~\ref{prelims}.

The tool referred to above is the double-bosonisation construction of Majid (\cite{BraidedLie}) for Lie bialgebras, which makes this work possible.  Assuming all our objects to be finite-dimensional, we take as input to the construction a quasitriangular Lie bialgebra $\smun{\Lie{g}}{0}$ and a braided-Lie bialgebra $\Lie{b}$.  One then obtains a new quasitriangular Lie bialgebra $\Lie{g} \iso \dbos{b}{\smun{g}{\mathrm{0}}}{\op{\dual{b}}}$.  (In general, one takes two dually paired braided-Lie bialgebras but we always take the usual dual in the finite-dimensional case.)  With this we can then ask two questions: \begin{enumerate} \item Given a Lie bialgebra $\Lie{g}$ and a sub-bialgebra $\smun{\Lie{g}}{0}$, \begin{enumerate} \item can we find a braided-Lie bialgebra $\Lie{b}$ such that we may reconstruct $\Lie{g}$ by double-bosonisation and \item if so, what is its structure? \end{enumerate} \item Given a Lie bialgebra $\smun{\Lie{g}}{0}$, for which choices of $\Lie{b}$ do we obtain something ``interesting'' on taking the double-bosonisation? \end{enumerate} 

\noindent For the first, we concentrate on the semisimple case with $\rank \smun{\Lie{g}}{0} = \rank \Lie{g} -1$.  We refer to this as the \emph{corank one} case.  In this situation, as shown by Majid, we have a positive answer to the first part.  Of course, we may repeat the process to answer the question for higher coranks.  The answer to the second part is the content of Section~\ref{deletion} and the Appendix of this paper.

An example of a negative answer to the first part is the \emph{corank zero} case: $\rank \smun{\Lie{g}}{0} = \rank \Lie{g}$.  In fact, the cases are quite different in character.  Although we may attack both representation-theoretically by restricting the adjoint representation of the larger algebra to the chosen subalgebra, when we consider the brackets on the modules we obtain, we have very different behaviour.  In the corank one case we have $\integ$-gradings, but in the corank zero case we see $\integ/(2)$- and $\integ/(3)$-gradings.  We cannot use double-bosonisation to reconstruct the larger algebra in the corank zero case, as it cannot reproduce these gradings by finite groups.  We note that the construction attributed to Freudenthal in Chapter 22 of \cite{FultonHarris} deals with the $\integ/(3)$-graded case purely on the level of Lie algebras.

For ``interesting'' in the second of the above questions, we take the (finite-dimensional) simple complex Lie algebras, with the Drinfel\cprime d--Sklyanin Lie cobracket.  We consider in Section~\ref{induction} the obvious induction scheme coming from Dynkin diagrams, namely, since we know that simple Lie algebras have connected Dynkin diagrams, adding a new node to a diagram in as many ways as is possible.  Clearly not all choices are allowed and we give some conditions and results, purely algebraic in nature, which control this.  We do not yet have a complete list of such conditions, so we cannot re-prove the classification of the simple complex Lie algebras, but we do have enough necessary conditions to exclude some possibilities immediately.  In particular, we discuss the obstructions to inducing in this way to the exceptional Lie algebras.

There has long been interest in the exceptional series of simple complex Lie algebras, for a variety of reasons in both mathematics and physics.  One goal has been to find a unified construction for these algebras.  Double-bosonisation goes beyond this, giving a unified construction of \emph{all} simple complex Lie algebras.  As we show here, it does so in a way that is completely compatible with the natural inclusions of Dynkin diagrams.  We refer to this scheme as \emph{Lie induction}.  Fully understanding Lie induction across the range of cases it encompasses is a programme that goes beyond this paper.  Here we set out some underlying theory and demonstrate the method for the exceptional series.

This, on its own, is a sufficient motivation to study double-bosonisation but there are others.  Firstly, it confirms the necessity of considering Lie bialgebras---not just Lie algebras---and their braided versions.  Also, through our calculations here, we have identified a number of new examples of braided-Lie bialgebra structures.

A second important goal has been to understand why we obtain only the infinite series and exceptionals we know in the classification.  The original proof of the classification is essentially a combinatorial analysis of a set of geometric objects, namely root systems.  Lie induction via double-bosonisation deals directly with the algebras and their representation theory.  In doing so, it emphasises the strong relationship between the structure theory and the representation theory, which has become a theme in algebra.  Using Lie induction, we are able to formulate in a new way the question of the obstruction to finding, for example, a finite-dimensional Lie algebra of type $E_{9}$.  With the analysis in this paper, we are now in a position to answer this question for $E_{9}$.

However, this work has implications outside of Lie algebra theory.  In recent years, the study of the quantized enveloping (Hopf) algebra $\qea{g}$ associated to a Lie algebra $\Lie{g}$ has become one of the most significant among algebraists, geometers and physicists alike.  These objects are still not fully understood, although much progress has of course been made.  We refer the reader to \cite{Jantzen} and \cite{FQGT} for introductions to this topic.

There is a strong relationship between what we refer to as the semi-classical theory, that of Lie bialgebras, and the quantum, in the form of $\qea{g}$.  Most notably from our perspective, double-bosonisation has also been defined for Hopf algebras, in work of Majid (\cite{DoubleBos}) and (independently) Sommerh{\"a}user (\cite{Sommerhauser}).  So, the methods we develop here should carry over to the quantum setting, providing some new insights into the structure of $\qea{g}$.  In place of braided-Lie bialgebras, we have braided groups: these are Hopf algebras in braided categories.  As yet, relatively few examples are known but it is clear that an analysis such as we have carried out for Lie bialgebras will provide a large class of examples, likely new ones.  This provides considerable further motivation for our work here.

We begin by recalling the structures we will need in Section~\ref{prelims}, namely those of a quasitriangular Lie bialgebra and a braided-Lie bialgebra.  We state the results of Majid (\cite{BraidedLie}) defining the double-bosonisation construction and its natural induced quasitriangular structure, when the input is quasitriangular.

The main body of this paper falls into two complementary but intimately connected halves, Section~\ref{deletion} discussing deletion and Section~\ref{induction} on induction.  We know from \cite{BraidedLie} that to each choice of simple root in a semisimple Lie algebra $\Lie{g}$ we have associated a braided-Lie bialgebra $\Lie{b}$ such that we can recover $\Lie{g}$ by double-bosonisation from $\smun{\Lie{g}}{0}$, where $\smun{\Lie{g}}{0}$ is a Lie subalgebra with Dynkin diagram that of $\Lie{g}$ with the node corresponding to the chosen simple root deleted.  We refer to calculating these $\Lie{b}$ as \emph{deletion}.

Then the aim of Section~\ref{deletion} is to further analyse the structure of this $\Lie{b}$ and provide some tools for calculating $\Lie{b}$ explicitly.  Critical to this are Lemma~\ref{zgrading}, where we observe that associated to each simple root is a decomposition of $\Lie{g}$ which defines a $\integ$-grading, and Proposition~\ref{irred}, where we cite a result of Azad, Barry and Seitz (\cite{ABS}) which tells us that the homogeneous components of this grading are irreducible modules for the zeroth part (except the zeroth part itself, which is not irreducible).  

We have calculated the braided-Lie bialgebra structures associated to simple Lie algebras $\Lie{g}$ in the case where we delete a simple root corresponding to an extremal node in the Dynkin diagram, so that the subalgebra $\smun{\Lie{g}}{0}$ is simple.  That is, for deletion to a simple Levi subalgebra.  We record the full results of our calculations in an Appendix and refer the reader to Section~\ref{summary} for a summary table (Table~\ref{deletionstable}) of the modules we find for each deletion.  We also consider the r\^{o}le of the graph automorphisms of the Dynkin diagrams, in Section~\ref{automorphisms}.

In Section~\ref{induction}, we take the opposite view and ask for necessary conditions on braided-Lie bialgebras $\Lie{b}$ in general to obtain simple Lie bialgebras via double-bosonisation.  Some properties are immediate from our analysis of deletion, for example, that $\Lie{b}$ must be ($\integ$-)graded.  We list these in Section~\ref{necccond}.  

From this list, we have two properties of particular importance, namely that the irreducible graded components of $\Lie{b}$ should have all their weight spaces one-dimensional (property~(\ref{1dwtspcond}) on page~\pageref{1dwtspcond}) and should have at most two dominant weights (property~(\ref{2domwtscond}$'$)).  We call a module with these two properties a \emph{defining} module.  In Section~\ref{defmods}, we record a result communicated to us by Y. Bazlov which classifies the defining modules for the simple complex Lie algebras.  

This forms the basis of our analysis in Section~\ref{indandexcep} of some of the obstructions to the existence of simple exceptional Lie algebras other than those already known.  Specifically, we see that there exist no non-trivial defining modules for $E_{8}$, and so we cannot produce a finite-dimensional algebra of type $E_{9}$.  Analysis of possible inductions from $A_{8}$ and $D_{8}$ only reinforce this.  Although we obtain some possible candidates from each, we have none that are consistent.  We make similar analyses for $F_{5}$ and $G_{3}$, concluding that we have no possible $F_{5}$ but we do obtain candidates for $G_{3}$.  In fact, we have countably many such candidates, which of course cannot be simple.  We cannot exclude these by our current list of necessary conditions, showing that these are not sufficient.

We conclude with some remarks on further extensions of this work.

%Acknowledgements

The author would like to thank Shahn Majid and Steve Donkin for all their comments and guidance.  Particular thanks go to Gerhard R\"{o}hrle for providing the reference \cite{ABS} and Yuri Bazlov for Theorem~\ref{defmodsthm}.

\section{Preliminaries}\label{prelims}

Throughout, unless otherwise stated, we work over the field of complex numbers $\complex$, although many of the definitions and some of the results can be extended to fields $k$ of characteristic not two.  We will make further comment on this later.  We assume that the reader is familiar with the basics of the theory of semisimple Lie algebras and root systems, as can be found in \cite{Serre}, \cite{Humphreys} or \cite{FultonHarris}, for example.  Our notations will typically follow those of the first two of these.  In particular, we will usually label the simple Lie algebras by the Cartan labelling ($A_{l}$, $B_{l}$, etc., where $l$ is the rank) and use the numbering of the simple roots as given on page~58 of \cite{Humphreys}.  We also assume knowledge of the highest weight theory of representations.

For an algebra-subalgebra pair $(\Lie{g},\smun{\Lie{g}}{0})$, $\smun{\Lie{g}}{0} \subset \Lie{g}$, with $\Lie{g}$ and $\smun{\Lie{g}}{0}$ semisimple, define the \emph{corank} of $\smun{\Lie{g}}{0}$ in $\Lie{g}$ to be $\corank(\Lie{g},\smun{\Lie{g}}{0})=\rank(\Lie{g})-\rank(\smun{\Lie{g}}{0})$.  In this paper, we will only concern ourselves with corank one pairs.

Our notation for highest-weight $\Lie{g}$-modules will be $V(\lambda; \Lie{g})$ for a highest weight $\lambda$, unless the algebra in question is clear from the context, when we write $V(\lambda)$.  We may use a notation for a specific realisation of the module, for example $S^{+}$ for a positive spin representation.  The notation $V$ (with no $\lambda$) is reserved for the appropriate natural representation, unless otherwise stated.

We use $\tau$ to mean the tensor product flip map, \eg \[ \tau:V\tensor W \to W\tensor V,\ \tau(v\tensor w)=w\tensor v\ \mbox{for all}\ v\in V,\ w\in W, \] on any appropriate pair of vector spaces.  The adjoint action of a Lie algebra $\Lie{g}$ on itself can be extended naturally to tensor products as follows.  For $x,y,z\in \Lie{g}$, \[ \ad{x}{y\tensor z}=\ad{x}{y}\tensor z +y\tensor \ad{x}{z}. \]  We will use this throughout without further comment.  We use the term ad-invariant in the obvious way.

We adopt the Sweedler notation for elements of tensor products.  That is, we use upper or lower parenthesized indices to indicate the placement in the tensor product, \eg $\sum a_{(1)} \tensor a_{(2)} \tensor a_{(3)} \in \Lie{g} \tensor \Lie{g} \tensor \Lie{g}$.  We will usually omit the summation sign.

The definition of a Lie bialgebra is originally due to Drinfel\cprime d (\cite{Drinfeld1}, \cite{Drinfeld2}).  The idea is the same as that for Hopf algebras, where we have two structures dual to each other, compatible in a natural way.  It is worth commenting that Lie bialgebras form a richer class than Lie algebras: the choice of the cobracket, the dual structure to the bracket, is not usually unique.

\begin{definition}[{\cite{Drinfeld1}}] A Lie bialgebra is $(\Lie{g},\Lbracket{\ }{\ },\delta)$ where 
\begin{enumerate} 
\item $(\Lie{g},\Lbracket{\ }{\ })$ is a Lie algebra, 
\item $(\Lie{g},\delta)$ is a Lie coalgebra, that is, $\delta:\Lie{g}\to \Lie{g}\tensor \Lie{g}$ satisfies 
	\[ \begin{array}{lll} \delta+\tau\circ \delta = 0 & & \mbox{\textup{(anticocommutativity)}} \\
			(\delta \tensor \id)\circ \delta +\mbox{\textup{cyclic}}=0 & & 
							\mbox{\textup{(co-Jacobi identity)}}
	\end{array} \]  (Here,\ \textup{``cyclic''}\ refers to cyclical rotations of the three tensor product factors in $\Lie{g}\tensor \Lie{g} \tensor \Lie{g}$.)
\item we have a cohomological compatibility condition: $\delta$ is a 1-cocycle in $Z_{\mathrm{ad}}^{1}(\Lie{g},\Lie{g}\tensor\Lie{g})$.  Explicitly, \[ \delta(\Lbracket{x}{y})=\ad{x}{\delta y}-\ad{y}{\delta x}. \] 
\end{enumerate}
\end{definition}

\noindent Examining this definition, we see that if $\Lie{g}$ is a finite-dimensional Lie bialgebra, then $(\dual{\Lie{g}},\dual{\delta},\dual{\Lbracket{\ }{\ }})$ is also a finite-dimensional Lie bialgebra.  Here, $\dual{\delta}$ and $\dual{\Lbracket{\ }{\ }}$ are the bracket and cobracket, respectively, given by dualisation.

In many of the natural cases we wish to consider, the cobracket $\delta$ arises as the coboundary of an element $r\in \Lie{g}\tensor\Lie{g}$.  (Explicitly, $\delta x=\ad{x}{r}$ for all $x\in \Lie{g}$.)  If $\delta$ further satisfies $(\id \tensor \delta)r=\Lbracket{\smun{r}{13}}{\smun{r}{12}}$, we say that $(\Lie{g},r)$ is a quasitriangular Lie bialgebra.  Here, we write $\smun{r}{12}=r^{(1)}\tensor r^{(2)} \tensor 1$, etc., with summation understood and the indices showing the placement in the triple tensor product $\Lie{g} \tensor \Lie{g} \tensor \Lie{g}$.  The bracket is taken in the common factor, so $\Lbracket{\smun{r}{13}}{\smun{r}{12}}=\Lbracket{r^{(1)}}{{r'}^{(1)}} \tensor {r'}^{(2)} \tensor {r}^{(2)}$ with $r'$ a second copy of $r$. 

To construct a quasitriangular Lie bialgebra, it is sufficient to find an element $r\in \Lie{g} \tensor \Lie{g}$ satisfying the classical Yang--Baxter equation and with ad-invariant symmetric part.  Then we take the coboundary $\partial r$ for $\delta$.  The classical Yang--Baxter equation, in the Lie setting, is \[ \Lbbracket{r}{r} \defeq \Lbracket{\smun{r}{12}}{\smun{r}{13}} + \Lbracket{\smun{r}{12}}{\smun{r}{23}} + \Lbracket{\smun{r}{13}}{\smun{r}{23}} = 0. \]  The bracket $\Lbbracket{\ }{\ }$ is the Schouten bracket, the natural extension of the bracket to these tensor spaces.

Considering the symmetric part of $r$, $2r_{+}\defeq r+\tau(r)$, we can distinguish two further cases.  Firstly, if $2r_{+}=0$ we say $(\Lie{g},r)$ is triangular.  Secondly, considering $2r_{+}$ as a map $\dual{\Lie{g}} \to \Lie{g}$, if this map is surjective, we say $(\Lie{g},r)$ is factorisable.  We refer the reader to the paper of Reshetikhin and Semenov-Tyan-Shanski\u{\i} (\cite{ReshSTS}).

We now consider the braided version of Lie bialgebras, as defined by Majid in \cite{BraidedLie}.  Here we consider the module category $\modcat{\Lie{g}}$ of a quasitriangular Lie bialgebra $\Lie{g}$ and objects in this category possessing a $\Lie{g}$-covariant Lie algebra structure.  Following the line suggested by the theory of braided groups, we associate to these objects a braiding-type map generalising the usual flip.  If $\Lie{b}$ is a $\Lie{g}$-covariant Lie algebra in the category $\modcat{\Lie{g}}$, we define the infinitesimal braiding of $\Lie{b}$ to be the operator $\psi:\Lie{b}\tensor \Lie{b} \to \Lie{b}\tensor\Lie{b}$, $\psi(a\tensor b)=2r_{+} \laction (a\tensor b-b\tensor a)$ where $\laction$ is the left action of $\Lie{g}$ on $\Lie{b}$ extended to the tensor products.  In fact, $\psi$ is a 2-cocycle in $Z_{\mathrm{ad}}^{2}(\Lie{b},\Lie{b}\tensor \Lie{b})$.

\begin{definition}[{\cite{BraidedLie}}]\label{blba} A braided-Lie bialgebra $(\Lie{b},\Lbracket{\ }{\ }_{\Lie{b}},\underline{\delta})$ is an object in $\modcat{\Lie{g}}$ satisfying the following conditions.  \begin{enumerate} \item $(\Lie{b},\Lbracket{\ }{\ }_{\Lie{b}})$ is a $\Lie{g}$-covariant Lie algebra in the category. \item $(\Lie{b},\underline{\delta})$ is a $\Lie{g}$-covariant Lie coalgebra in the category. \item $\dform \underline{\delta}=\psi$. \end{enumerate} \end{definition}

\noindent We can now state the theorem which provides the construction we use in this paper.  Let $\Lie{g}$ be a quasitriangular Lie bialgebra.

\begin{theorem}[{\cite[Theorem 3.10]{BraidedLie}}]\label{dbos} For dually paired braided-Lie bialgebras $\Lie{b},\Lie{c}\in \modcat{\Lie{g}}$ the vector space $\Lie{b} \dsum \Lie{g} \dsum \Lie{c}$ has a unique Lie bialgebra structure $\dbos{b}{g}{\op{c}}$, the double-bosonisation, such that $\Lie{g}$ is a sub-Lie bialgebra, $\Lie{b},\Lie{\op{c}}$ are Lie subalgebras, and \[ \begin{array}{c} \Lbracket{\xi}{x}=\xi \laction x, \quad \Lbracket{\xi}{\phi}=\xi \laction \phi \\ \\  \Lbracket{x}{\phi} = x_{\underline{(1)}}\ip{\phi}{x_{\underline{(2)}}}+\phi_{\underline{(1)}}\ip{\phi_{\underline{(2)}}}{x}+2r_{+}^{(1)}\ip{\phi}{r_{+}^{(2)}\laction x} \\ \\ \delta x=\underline{\delta}x+r^{(2)}\tensor r^{(1)}\laction x-r^{(1)}\laction x\tensor r^{(2)} \\ \\ \delta \phi=\underline{\delta}\phi + r^{(2)}\laction \phi \tensor r^{(1)} - r^{(1)}\tensor r^{(2)}\laction \phi \end{array} \] \noindent for all $x\in \Lie{b}$, $\xi\in \Lie{g}$ and $\phi \in \Lie{c}$.  Here $\underline{\delta}x=x_{\underline{(1)}}\tensor x_{\underline{(2)}}$. \end{theorem}

\noindent Moreover, the double-bosonisation is always quasitriangular when we take $\Lie{c}=\dual{b}$, as we see from the following proposition.

\begin{proposition}[{\cite[Proposition 3.11]{BraidedLie}}]\label{doubleqtstr} Let $\Lie{b}\in \modcat{\Lie{g}}$ be a finite-dimensional braided-Lie bialgebra with dual $\dual{\Lie{b}}\!$.  Then the double-bosonisation $\dbos{b}{g}{\op{\dual{b}}}$ is quasitriangular with \[ r^{\mathrm{new}}=r+\sum_{a} f^a \tensor e_a \] where $\{ e_a \}$ is a basis of $\Lie{b}$ and $\{ f^a \}$ is a dual basis, and $r$ is the quasitriangular structure of $\Lie{g}$.  If $\Lie{g}$ is factorisable then so is the double-bosonisation.
\end{proposition}

In this paper, we concentrate on the Lie algebra structure rather than the coalgebra structure.  For each semisimple Lie algebra, there exists a canonical quasitriangular structure, the Drinfel\cprime d--Sklyanin solution (see \cite[Chapter 22]{QGP} for more details) and, unless otherwise stated, we use this choice.

\section{Deletion}\label{deletion}

Our aim is to associate to each corank-one pair of finite-dimensional semisimple complex Lie bialgebras $(\Lie{g},\smun{\Lie{g}}{0})$, $\smun{\Lie{g}}{0} \subset \Lie{g}$, a $\smun{\Lie{g}}{0}$-module $\Lie{b}=\Lie{b}(\Lie{g},\smun{\Lie{g}}{0})$ which, with the additional structure of a braided-Lie bialgebra, realises the induction from $\smun{\Lie{g}}{0}$ to $\Lie{g}$ given by an isomorphism of Lie bialgebras of $\Lie{g}$ with $\dbos{b}{\widetilde{\smun{g}{\mathrm{0}}}}{\op{\dual{b}}}$.  Here, $\widetilde{\smun{\Lie{g}}{0}}$ denotes a suitable central extension of $\smun{\Lie{g}}{0}$ which raises the rank by one.  

To do this, we use a combination of structure theory and representation theory to give some general tools, described below, and we give our explicit calculations in an Appendix.  It is clear that without loss of generality we may assume the larger algebra $\Lie{g}$ is simple.  However, we do not assume that the subalgebra $\smun{\Lie{g}}{0}$ is simple, unless otherwise stated.

\subsection{Gradings associated to simple roots}\label{grading}

We exhibit here a $\integ$-grading associated to each choice of simple root in a Lie algebra $\Lie{g}$.  It is this grading that will give us most of the information we need to determine the braided-Lie bialgebra $\Lie{b}$ discussed above.

Choose a Cartan subalgebra $\Lie{h}$ for $\Lie{g}$, a simple complex Lie algebra, and let $R$ be the associated root system.  Let $S=\{ \alpha_{1}, \ldots , \alpha_{l} \}$ be a base of simple positive roots for $R$ where $l=\dim \Lie{h} = \rank(\Lie{g})$.  Choose a Weyl basis for $\Lie{g}$, as follows: $\Lie{g}$ is generated by elements $H_{i} \in \Lie{h}$ corresponding to the $\alpha_{i}$ and elements $X_{i}^{+} \in \Lie{g}^{\alpha_{i}}$, $X_{i}^{-} \in \Lie{g}^{-\alpha_{i}}$ satisfying $\Lbracket{X_{i}^{+}}{X_{i}^{-}}=H_{i}$.  In particular, we have the Weyl relations \begin{align*} & \Lbracket{H_{i}}{X_{j}^{+}}=C_{ij}X_{j}^{+}, \\ & \Lbracket{H_{i}}{X_{j}^{-}}= -C_{ij}X_{j}^{-} \\ & \Lbracket{X_{i}^{+}}{X_{j}^{-}} = 0\quad \text{if $i\neq j$} \end{align*} where $C$ is the Cartan matrix associated to $\Lie{g}$.  The full basis is \[ \{ H_{i}, X_{i}^{+}, X_{i}^{-} \mid 1 \leq i \leq l \} \union \{ X_{\alpha}^{+}, X_{\alpha}^{-} \mid X_{\alpha}^{\pm} \in \Lie{g}^{\pm \alpha}, \alpha \in R^{+}\setminus S \}, \] where $R^{+}$ is the set of positive roots in $R$.
				   
We will want to consider subsets of the negative roots to define $\Lie{b}$ and we will work with the coordinate system provided by $S$, \ie if $\alpha$ is a root we can write $\alpha = \sum_{i=1}^{l} k_{i}\alpha_{i}$ and we have all $k_i \geq 0$ if and only if $\alpha \in R^{+}$ and conversely all $k_i \leq 0$ if and only if $\alpha \in R^{-}$, the set of negative roots.  Define the support of $\alpha \in R$ to be $\supp{\alpha} = \{ \alpha_{i} \in S \mid k_i \neq 0 \}$ when $\alpha$ is written in this way.  We will call $k_i$ the multiplicity of $\alpha_{i}$ in $\alpha$: $\mult{i}{\alpha} \defeq k_{i}$.  Finally, define the height of $\alpha$ to be $\height{\alpha}=\sum_{i=1}^{l} k_{i}$, if $\alpha = \sum_{i=1}^{l} k_{i}\alpha_{i}$.  In particular, the simple roots $\alpha_{i} \in S$ have height one.
				   
Since we have assumed $\Lie{g}$ to be simple, there exists a unique root $\Lambda$ in $R^{+}$ with maximal height, \ie $\height{\alpha} < \height{\Lambda}$ for all $\alpha \neq \Lambda$, $\alpha \in R$ (see, for example, \cite[Lemma 10.4.A]{Humphreys}).  We will call $\Lambda$ the highest root in $R$ (or $\Lie{g}$).  We recall that $\Lambda$ is also the highest weight vector in the adjoint representation.  The coordinate expression as a root for $\Lambda$, $\Lambda = \sum_{i=1}^{l} m_i \alpha_{i}$, may therefore be obtained from the expression for $\Lambda$ in the basis of fundamental weights via multiplication by $\inv{C}$.  For later use, we record these dual expressions for the irreducible root systems (labelled by the Cartan type) in Table~\ref{highroots}.  In what follows, we use parentheses $(\, .\, , \ldots , .\, )$ for vectors in the root basis provided by $S$ and square brackets $[\, .\, , \ldots , .\, ]$ for weights using the fundamental weights $\{ \omega_{i} \mid 1 \leq i \leq l \}$ (dual to $S$) as basis.  We will use the notation $\omega_{0}$ for the zero weight $[0, \ldots , 0]$.

\begin{table}
\begin{tabular}{l|l|l|l@{$\ =\ $}l}
Type & Highest root, $\Lambda$ & $\height{\Lambda}$ & \multicolumn{2}{l}{\parbox[t]{5cm}{$\omega_{\mathrm{ad}}$, highest weight of adjoint representation}} \\ \hline
$A_l$ & $(1,1,\ldots,1)$       & $l$    & $[1,0,0,\ldots,0,1]$ & $\omega_1+\omega_l$ \\
$B_l$ & $(1,2,2,\ldots,2)$     & $2l-1$ & $[0,1,0,\ldots,0]$   & $\omega_{2}$        \\
$C_l$ & $(2,2,\ldots,2,1)$     & $2l-1$ & $[2,0,0,\ldots,0]$   & $2\omega_{1}$       \\
$D_l$ & $(1,2,2,\ldots,2,1,1)$ & $2l-3$ & $[0,1,0,\ldots,0]$   & $\omega_{2}$        \\
$E_6$ & $(1,2,2,3,2,1)$        & $11$   & $[0,1,0,0,0,0]$      & $\omega_{2}$        \\
$E_7$ & $(2,2,3,4,3,2,1)$      & $17$   & $[1,0,0,0,0,0,0]$    & $\omega_{1}$        \\
$E_8$ & $(2,3,4,6,5,4,3,2)$    & $29$   & $[0,0,0,0,0,0,0,1]$  & $\omega_{8}$        \\
$F_4$ & $(2,3,4,2)$            & $11$   & $[1,0,0,0]$          & $\omega_{1}$        \\
$G_2$ & $(3,2)$                & $5$    & $[0,1]$              & $\omega_{2}$
\end{tabular}
\caption{Expressions for highest roots in the irreducible root systems\label{highroots}}
\end{table}

Observe also that $\height{\Lambda}=h-1$, where $h$ is the Coxeter number of $\Lie{g}$ (\cite[Ch. 6, Prop. 1.11.31]{Bourbaki}).

Let $J$ be a subset of $\{1, \ldots , l\}$.  The root deletion of $J$ is the 4-tuple $(\Lie{g},J,\smun{\Lie{g}}{0},\iota)$, where $\smun{\Lie{g}}{0}$ is the subalgebra of $\Lie{g}$ generated by the $3(l-|J|)$ generators $\{ H_{i}, X_{i}^{+}, X_{i}^{-} \mid i \not\in J \}$ and $\iota:\smun{\Lie{g}}{0} \inj \Lie{g}$ is the embedding of $\smun{\Lie{g}}{0}$ in $\Lie{g}$ defined by this choice of generators for $\smun{\Lie{g}}{0}$.  In the case when $|J|=1$, $J=\{ \alpha_{d} \}$ we write $(\Lie{g},d,\smun{\Lie{g}}{0},\iota)$.  Clearly, the Dynkin diagram for $\smun{\Lie{g}}{0}$ is given by deleting the nodes in the Dynkin diagram for $\Lie{g}$ corresponding to the $\alpha_{j}$, $j\in J$.  The map $\iota$ defines an embedding of the Dynkin diagram of $\smun{\Lie{g}}{0}$ into that for $\Lie{g}$ in the obvious way.

We now restrict to the case $|J|=1$, $J=\{ \alpha_{d} \}$, \ie the deletion of one simple root.  Let $\Lie{g}$ be a finite-dimensional complex simple Lie algebra.  

\begin{lemma}\label{zgrading} Associated to each simple root $\alpha_{d} \in S$, there is a $\integ$-grading of $\Lie{g}$ given by the $\alpha_{d}$-multiplicity as follows. Define $\mult{d}{X^{\pm}_{\alpha}} = \mult{d}{\alpha}$, $\alpha \in R$, and $\mult{d}{H_{i}} = 0$ for all $i=1,\ldots, l$.  Set \[ \smun{\Lie{g}}{[i]}=\setspan{\complex}{x\in \Lie{g} \mid \mult{d}{x}=i}, \] with the convention $\setspan{\complex}{\emptyset}=\{0\}$.  Then $\Lie{g}=\bigdsum_{i\in \integ} \smun{\Lie{g}}{[i]}$. \end{lemma}

\begin{proof} This is immediate from the additivity of $\mult{d}{-}$, coming from the additivity in the root system. \end{proof}

\noindent Note that this is not the usual $\integ$-grading of a finite-dimensional simple Lie algebra, with $\Lie{g}$ as the zero part and all other components zero.  In the above grading, the zero part is $\smun{\Lie{g}}{[0]}=\widetilde{\Lie{g}_{d}}$, a central extension of the subalgebra\footnote{In the above notation, this subalgebra is $\smun{\Lie{g}}{0}$ but we now denote it by $\Lie{g}_{d}$ to avoid confusion with the zero part of the grading and to emphasise its dependence on the deletion.} $\Lie{g}_{d} \subset \Lie{g}$ generated by all the generators of $\Lie{g}$ except $H_{d}$, $X^{+}_{d}$ and $X^{-}_{d}$.  The number of non-zero graded components is $2\cdot\mult{d}{\Lambda}+1$ ($\Lambda$ the highest root in $\Lie{g}$) and we see from Table~\ref{highroots} that we have $1\leq \mult{d}{\Lambda}=m_{d} \leq 6$ in general and $m_{d}\leq 3$ if $d$ is chosen such that $\Lie{g}_{d}$ is simple.

The most important property of this grading is that it gives the branching (\ie restriction) of the adjoint representation of $\Lie{g}$ to the subalgebra $\Lie{g}_{d}$.

\begin{proposition}\label{irred} For $i\neq 0$, $\smun{\Lie{g}}{[i]}$ is an irreducible $\Lie{g}_{d}$-module; $\smun{\Lie{g}}{[0]}=\Lie{g}_{d} \dsum \complex$ as $\Lie{g}_{d}$-modules. \end{proposition}

\begin{proof} The action of $\Lie{g}_{d}$ is induced by the bracket in $\Lie{g}$ and it is then clear that the $\smun{\Lie{g}}{[i]}$ are $\Lie{g}_{d}$-modules by the grading property.

That the $\smun{\Lie{g}}{[i]}$, $i\neq 0$, are irreducible may be deduced from a result of Azad, Barry and Seitz (\cite{ABS}).  Their results concern algebraic groups over more general fields but the parts we need are root system-theoretic and so carry across immediately.  The appropriate theorem in their paper is Theorem 2. \end{proof}

\begin{jegnote}{} We may observe that for $i=\pm 1, \pm m_{d}$, the irreducibility of $\smun{\Lie{g}}{[i]}$ is immediate.  The modules $\smun{\Lie{g}}{[\pm 1]}$ have a primitive generator, namely $X^{\pm}_{d}$; for $\smun{\Lie{g}}{[-1]}$, this highest weight vector has highest weight given by the negative of the $d$\/th row of the Cartan matrix for $\Lie{g}$ with the $d$\/th column deleted and re-ordered according to that induced by the embedding $\iota:\smun{\Lie{g}}{0} \inj \Lie{g}$.   The modules $\smun{\Lie{g}}{[\pm m_{d}]}$ have a unique lowest weight vector, $X^{\pm}_{\Lambda}$.  These observations are useful for the calculations we perform later. \end{jegnote}

\vspace{1em}
The above grading is related to double-bosonisation as follows.  Let $\Lie{n}^{-}$ be the standard negative Borel subalgebra of $\Lie{g}$, so $\Lie{n}^{-} = \Lie{h} \dsum \sum_{\alpha \in R^{-}} \Lie{g}^{\alpha}$.  Let $\Lie{b}$ be the Lie ideal of $\Lie{n}^{-}$ generated by $X_{d}^{-}$.  A basis for $\Lie{b}$ is $\{ X_{\alpha}^{-} \mid \alpha_{d} \in \supp{\alpha} \}$ and the subalgebra $\Lie{f}$ of $\Lie{n}^{-}$ generated by the set $\{ X_{i}^{-} \mid i \neq d \}$ has basis $\{ X_{\alpha}^{-} \mid \alpha_{d} \not\in \supp{\alpha} \}$.  Then we have the following.

\begin{proposition}\label{gradedbraided} Let $\Lie{g}$ be a finite-dimensional simple complex Lie bialgebra.  Choose a simple root of $\Lie{g}$, $\alpha_{d}$.  Then we have the decomposition  \[ \dbos{b}{\widetilde{{g}_{\mathit{d}}}}{\op{\dual{b}}} \] with $\Lie{g}_{d}$ generated by all the generators of $\Lie{g}$ except $H_{d}$, $X^{+}_{d}$ and $X^{-}_{d}$ and $\Lie{b}=\bigdsum_{i<0} \smun{\Lie{g}}{[i]}$, a $\integ$-graded braided-Lie bialgebra.  
\end{proposition}

\begin{proof} This follows from Proposition~4.5 of \cite{BraidedLie} and the definition of the grading associated to $\alpha_{d}$ in Lemma~\ref{zgrading}.
\end{proof}

\subsection{Automorphisms}\label{automorphisms}

Note that in both deletion and induction, we need to take account of the existence of graph automorphisms of some of the Dynkin diagrams associated to the simple Lie algebras, which we will call diagram automorphisms.  In deletion, we see certain symmetries appearing in the results of our calculations.  There are relatively few automorphisms to take care of---the list of simple Lie algebras with non-trivial automorphism group is as follows: $A_{l}$ (with automorphism group $\mathrm{Aut}\ \Lie{g}=\integ/(2)$), $D_{4}$ ($S_3$), $D_{l}$, $l\geq 5$ ($\integ/(2)$) and $E_{6}$ ($\integ/(2)$).  Observe that these are all simply-laced algebras, that is, there is only one root length in the root system.

As a result of the existence of these automorphisms, we want to consider certain deletions (as defined in Section~\ref{grading}) equivalent, denoted $\simeq$.  It is diagram automorphisms that lead us to insist on specifying the embedding $\iota$ as part of the deletion data but we now record which give essentially the same data.  By ``essentially'', we mean that we may not find the same modules but may find their duals (where this is different).  There will be $| \mathrm{Aut}\ \Lie{g} |\cdot | \mathrm{Aut}\ \smun{\Lie{g}}{0} |$ equivalent deletions $(\Lie{g},-,\smun{\Lie{g}}{0},-)$.  We have 
\begin{enumerate}
\renewcommand{\labelenumi}{\roman{enumi})}
\item $\Lie{g}=A_{l+1}$, $\smun{\Lie{g}}{0}=A_{l}$: \begin{align*} (A_{l+1},1,A_{l},i\mapsto i+1) & \simeq (A_{l+1},1,A_{l},i \mapsto l-i+2) \\ & \simeq (A_{l+1},l,A_{l},\id) \\ & \simeq (A_{l+1},l,A_{l},i\mapsto l-i+1) \end{align*}
\item $\Lie{g}=D_{4}$, $\smun{\Lie{g}}{0}=A_{3}$: \begin{align*} (D_4,1,A_3,\left(\begin{array}{lll} \Ss 1 & \Ss 2 & \Ss 3 \\ \Ss 3 & \Ss 2 & \Ss 4 \end{array} \right)) & \simeq (D_4,1,A_3,\left( \begin{array}{lll} \Ss 1 & \Ss 2 & \Ss 3 \\ \Ss 4 & \Ss 2 & \Ss 3 \end{array} \right)) \\ & \simeq (D_4,3,A_3,\left( \begin{array}{lll} \Ss 1 & \Ss 2 & \Ss 3 \\ \Ss 1 & \Ss 2 & \Ss 4 \end{array} \right)) \\ & \simeq (D_4,3,A_3,\left( \begin{array}{lll} \Ss 1 & \Ss 2 & \Ss 3 \\ \Ss 4 & \Ss 2 & \Ss 1 \end{array} \right)) \\ & \simeq (D_4,4,A_3,\id) \\ & \simeq (D_4,4,A_3,\left( \begin{array}{lll} \Ss 1 & \Ss 2 & \Ss 3 \\ \Ss 3 & \Ss 2 & \Ss 1 \end{array} \right)) \end{align*}
\item $\Lie{g}=D_{l+1}$, $\smun{\Lie{g}}{0}=A_{l}$: \begin{align*} (D_{l+1},l-1,A_{l},&\left( \begin{array}{llllll} \Ss 1 & \Ss 2 & \Ss 3 & \Ss \cdots & \Ss l-1 & \Ss l \\ \Ss 1 & \Ss 2 & \Ss 3 & \Ss \cdots & \Ss l-1 & \Ss l+1 \end{array} \right)) \\ & \simeq (D_{l+1},l-1,A_{l},\left( \begin{array}{llllll} \Ss 1 & \Ss 2 & \Ss 3 & \Ss \cdots & \Ss l-1 & \Ss l \\ \Ss l+1 & \Ss l-1 & \Ss l-2 & \Ss \cdots & \Ss 2 & \Ss 1 \end{array} \right)) \\ & \simeq (D_{l+1},l,A_{l},\id) \\ & \simeq (D_{l+1},l,A_{l},i\mapsto l-i+1) \end{align*}
\item $\Lie{g}=E_{6}$, $\smun{\Lie{g}}{0}=D_{5}$: \begin{align*}  (E_6,1,D_5,\left( \begin{array}{lllll} \Ss 1 & \Ss 2 & \Ss 3 & \Ss 4 & \Ss 5 \\ \Ss 6 & \Ss 5 & \Ss 4 & \Ss 3 & \Ss 2 \end{array} \right)) & \simeq (E_6,1,D_5,\left( \begin{array}{lllll} \Ss 1 & \Ss 2 & \Ss 3 & \Ss 4 & \Ss 5 \\ \Ss 6 & \Ss 5 & \Ss 4 & \Ss 2 & \Ss 3 \end{array} \right)) \\ & \simeq (E_6,6,D_5,\left( \begin{array}{lllll} \Ss 1 & \Ss 2 & \Ss 3 & \Ss 4 & \Ss 5 \\ \Ss 1 & \Ss 3 & \Ss 4 & \Ss 2 & \Ss 5 \end{array} \right)) \\ & \simeq (E_6,6,D_5,\left( \begin{array}{lllll} \Ss 1 & \Ss 2 & \Ss 3 & \Ss 4 & \Ss 5 \\ \Ss 1 & \Ss 3 & \Ss 4 & \Ss 5 & \Ss 2 \end{array} \right)) \end{align*}
\end{enumerate}

\noindent We refer the reader to the Appendix (p.\:\pageref{appendix}) for a full explanation of our notation.

We note that the potentially interesting case of the triple symmetry in the diagram for $D_{4}$ does \emph{not} yield three different representations but in fact $\Lie{g}_{[-1]}=V(\omega_{2}; A_{3})=\Wedge{2}{(V)}$ in all cases.

\subsection{Summary of deletions}\label{summary}

We give here a summary of our calculations, with the details reserved for the Appendix.  Firstly, in the above we did not consider how we obtained $A_{1}=\Lie{sl}_{2}$, since $A_1$ does not have a simple semisimple subalgebra.  However, Majid observed in \cite{BraidedLie} that the procedure of deleting all the roots from a Lie algebra, leaving just the Cartan subalgebra, and its corresponding induction make sense and he gives general formul{\ae} there.  For completeness, we record this deletion for $A_1$ here.

\begin{description} 
\item[($A_{1}\complex$)] Deletion $(A_{1}, 1, \Lie{h}=\complex\cdot H, - )$

$\Lie{b}_{-1}$ is spanned by $X^{-}$ and we have as action $H\laction X^{-} = -2X^{-}$.  For the dual, $\Lie{b}_{1}$, we choose as basis $X^{+}$ with $X^{+}(X^{-})=-1$ (the negative of the usual choice).  Then $H\laction X^{+}=2X^{+}$.  We consider $\Lie{b}_{-1}=\complex \cdot X^{-}$ as a braided-Lie bialgebra with the zero bracket and cobracket and this induces the same for $\Lie{b}_{1}$.  Note that the infinitesimal braiding is also zero.

One may check that $\complex \cdot H$ with the zero bracket, quasitriangular structure $r=\frac{1}{4} H\tensor H$ and the above action gives the double-bosonisation ${\complex \cdot X^{-}} \rtimesdot {\complex \cdot H} \ltimesdot {\complex \cdot X^{+}} \iso \Lie{sl}_{2}=A_{1}$ with the Drinfel\cprime d--Sklyanin quasitriangular structure.  Here we do not need to make a further central extension.
\end{description}

We now give a table (Table~\ref{deletionstable}) summarising the remainder of our calculations, that is, for the (equivalence classes of) deletions $(\Lie{g},d,\smun{\Lie{g}}{0},\iota)$ with $\corank (\Lie{g},\smun{\Lie{g}}{0})=1$ and \emph{both} $\Lie{g}$ and $\smun{\Lie{g}}{0}$ simple.  From Proposition~\ref{gradedbraided}, we know that $\Lie{b}$, the braided-Lie bialgebra arising from deletion, is graded with irreducible components and it is these modules occurring in $\Lie{b}$ that we give here.  More details of the full braided-Lie bialgebra structure are given in the Appendix.  We also indicate the type of representation, \ie trivial, the natural representation, a spin representation, etc.

\pagebreak

\begin{landscape}

\begin{table}
\caption{Summary of deletions\label{deletionstable}}
\begin{tabular}{llll|l|llllll}
$\Lie{g}$ & $d$ & $\smun{\Lie{g}}{0}$ & $\iota$ & $m_{d}$ & \multicolumn{2}{l}{$\Lie{b}_{-1}=\Lie{g}_{[-1]}$} & \multicolumn{2}{l}{$\Lie{b}_{-2}=\Lie{g}_{[-2]}$} & \multicolumn{2}{l}{$\Lie{b}_{-3}=\Lie{g}_{[-3]}$} \\ \hline

$A_{l+1}$ & $l$ & $A_{l}$ & $\id$ & 1 & $\omega_l$ & natural & & & & \\

$B_{l+1}$ & 1 & $B_{l}$ & $i\mapsto i+1$ & 1 & $\omega_1$ & natural & & & & \\

$C_{l+1}$ & 1 & $C_{l}$ & $i\mapsto i+1$ & 2 & $\omega_1$ & natural & $\omega_0$ & trivial & & \\

$D_{l+1}$ & 1 & $D_{l}$ & $i\mapsto i+1$ & 1 & $\omega_1$ & natural & & & & \\

$E_{7}$ & 7 & $E_{6}$ & $\id$ & 1 & $\omega_6$ & & & & & \\

$E_{8}$ & 8 & $E_{7}$ & $\id$ & 2 & $\omega_7$ & & $\omega_0$ & trivial & & \\

$B_{l+1}$ & $l+1$ & $A_{l}$ & $i\mapsto l-i+1$ & 2 & $\omega_1$ & natural & $\omega_2$ & $\Wedge{2}{(\text{natural})}$ & & \\

$C_{l+1}$ & $l+1$ & $A_{l}$ & $i\mapsto l-i+1$ & 1 & $2\omega_1$ & $\Sym{2}{\text{natural}}$ & & & & \\

$D_{l+1}$ & $l+1$ & $A_{l}$ & $i\mapsto l-i+1$ & 1 & $\omega_2$ & $\Wedge{2}{(\text{natural})}$ & & & & \\

$E_{6}$ & 2 & $A_{5}$ & $\left( \begin{array}{lllll} \Ss 1 & \Ss 2 & \Ss 3 & \Ss 4 & \Ss 5 \\ \Ss 1 & \Ss 3 & \Ss 4 & \Ss 5 & \Ss 6 \end{array} \right)$ & 2 & $\omega_3$ & $\Wedge{3}{(\text{natural})}$ & $\omega_0$ & trivial & & \\

$E_{7}$ & 2 & $A_{6}$ & $\left( \begin{array}{llllll} \Ss 1 & \Ss 2 & \Ss 3 & \Ss 4 & \Ss 5 & \Ss 6 \\ \Ss 1 & \Ss 3 & \Ss 4 & \Ss 5 & \Ss 6 & \Ss 7 \end{array} \right)$ & 2 & $\omega_3$ & $\Wedge{3}{(\text{natural})}$ & $\omega_6$ & $\Wedge{6}{(\text{natural})}$ & & \\

$E_{8}$ & 2 & $A_{7}$ & $\left( \begin{array}{lllllll} \Ss 1 & \Ss 2 & \Ss 3 & \Ss 4 & \Ss 5 & \Ss 6 & \Ss 7 \\ \Ss 1 & \Ss 3 & \Ss 4 & \Ss 5 & \Ss 6 & \Ss 7 & \Ss 8 \end{array} \right)$ & 3 & $\omega_3$ & $\Wedge{3}{(\text{natural})}$ & $\omega_6$ & $\Wedge{6}{(\text{natural})}$ & $\omega_1$ & natural \\

$G_{2}$ & 1 & $A_{1}$ & $\left( \begin{array}{l} \Ss 1 \\ \Ss 2 \end{array} \right)$ & 3 & $\omega_1$ & natural & $\omega_0$ & trivial & $\omega_1$ & natural \\

$G_{2}$ & 2 & $A_{1}$ & $\id$ & 2 & $3\omega_1$ & $\Sym{3}{\text{natural}}$ & $\omega_0$ & trivial & & \\

$F_{4}$ & 1 & $C_3$ & $i \mapsto 5-i$ & 2 & $\omega_3$ & & $\omega_{0}$ & trivial & & \\

$F_{4}$ & 4 & $B_3$ & $\id$ & 2 & $\omega_3$ & spin & $\omega_1$ & natural & & \\

$E_{6}$ & 1 & $D_{5}$ & $i \mapsto 7-i$ & 1 & $\omega_4$ & positive spin & & & & \\

$E_{7}$ & 1 & $D_{6}$ & $i \mapsto 8-i$ & 2 & $\omega_5$ & negative spin & $\omega_0$ & trivial & & \\

$E_{8}$ & 1 & $D_{7}$ & $i \mapsto 9-i$ & 2 & $\omega_6$ & positive spin & $\omega_1$ & natural & &
\end{tabular}

\end{table}
\end{landscape}

\pagebreak

\section{Induction}\label{induction}

We now begin the programme to analyse the classification of the simple Lie algebras using the representation-theoretic approach of Lie induction.  This gives a somewhat different perspective to the usual method for classifying the simple Lie algebras using the geometry and combinatorics of the root systems.  

Our first task is to see to what general principles we can extract from the above analysis of deletions to give necessary conditions for braided-Lie bialgebras to induce simple Lie algebras.  More precisely, we wish to analyse the properties required by a braided-Lie bialgebra $\Lie{b}$ in the module category of a finite-dimensional simple Lie algebra $\smun{\Lie{g}}{0}$ so that the double-bosonisation $\dbos{b}{\widetilde{\smun{g}{\mathrm{0}}}}{\op{\dual{b}}}$ is again finite-dimensional and simple.

In particular, we would like to understand what the obstructions are that limit the classification to the known series and exceptionals.  This could suggest whether or not relaxing certain axioms would alter the classification, \eg using quasi-Lie algebras.  To do this, we use a classification of irreducible modules satisfying two key necessary conditions.  This determines the modules which may appear as irreducible components in $\Lie{b}$, which we know to be graded.  Then if no such modules exist for a given simple $\smun{\Lie{g}}{0}$, there can be no induction.

\subsection{Necessary conditions on $\Lie{b}$}\label{necccond}

The key idea is that we are considering modules which are potential subsets of roots in irreducible root systems and the following conditions come from this and the structure discussed in the previous section.  Firstly, we recall (Lemma~\ref{zgrading}) that we can consider a simple Lie algebra $\Lie{g}$ to be $\integ$-graded by choosing a simple root $\alpha_{d}$ and grading by $\mult{d}{X^{\pm}_{\alpha}} = \mult{d}{\alpha}$, $\mult{d}{H_{i}}=0$.  So we have the condition \begin{property} \Lie{b}\ \text{should be a finite-dimensional graded Lie algebra.} \end{property}  That is, $\Lie{b}=\bigdsum_{j=-m}^{-1} \Lie{b}_{j}$, $\Lbracket{\Lie{b}_{j}}{\Lie{b}_{k}} \subseteq \Lie{b}_{j+k}$ (possibly zero) with $m<\infty$.  Next, the homogeneous parts should be irreducible: \begin{property} \Lie{b}_{j}\ \text{should be irreducible for all}\ j=-1, \ldots, -m. \end{property}  This comes from the theorems of \cite{ABS}.  For conditions on the candidates for the $\Lie{b}_{j}$, we look to the underlying irreducible root system.  Any root system of a Lie algebra has one-dimensional root spaces so we require \begin{property}\label{1dwtspcond} \Lie{b}_{j}\ \text{has all its weight spaces one-dimensional.} \end{property}  Also, there can be at most two root lengths and the roots of the same length and height must be conjugate under the Weyl group of $\smun{\Lie{g}}{0}$ (\cite{ABS}) so \begin{property} \label{2domwtscond} \Lie{b}_{j}\ \text{has at most two Weyl group orbits.} \end{property}  This can be rephrased in terms of dominant weights: \begin{equivproperty} \Lie{b}_{j}\ \text{has at most two dominant weights.} \end{equivproperty}  We will see that the conditions (\ref{1dwtspcond}) and (\ref{2domwtscond}) are very restrictive when combined.  We say a module is \emph{defining} if it satisfies (\ref{1dwtspcond}) and (\ref{2domwtscond}) (equivalently, (\ref{1dwtspcond}) and (\ref{2domwtscond}$'$)).

We also need a property related to the grading on $\Lie{b}$.  When discussing calculating deletions, we observe that $\Lie{b}_{-2}$ must be a submodule of $\Wedge{2}{(\Lie{b}_{-1})}$, by considering the module map $\Lbracket{\ }{\ }|_{\Lie{b}_{-1}}$, and similarly for higher indices.  So, we require \begin{property}\label{suboftensorprodcond} \Lie{b}_{k}\ \text{occurs as a submodule of}\ \Lie{b}_{i} \tensor \Lie{b}_{j}\ \text{for all}\ i,j\ \text{such that}\ i+j=k. \end{property}  This follows from considering the bracket maps, which will be module maps, and Schur's lemma.  If $i=j=k/2$, we require that $\Lie{b}_{k}$ occurs as a submodule of $\Wedge{2}{(\Lie{b}_{i})}$.

We will classify the defining modules for the simple Lie algebras in the next section but we can immediately see that the trivial module satisfies conditions (\ref{1dwtspcond}) and (\ref{2domwtscond}) and so is a candidate.  However, the following theorem discounts this possibility.
 
\begin{theorem}\label{nottrivial} Let $\Lie{g}$ be a finite-dimensional simple quasitriangular complex Lie bialgebra and $\complex$ be its trivial representation.  Then \[ \dbos{\complex}{\widetilde{g}}{\op{\dual{\complex}}} \iso \Lie{g}\dsum \Lie{sl}_{2}(\complex) \] as Lie bialgebras, where $\Lie{sl}_{2}(\complex)$ has the Drinfel\cprime d--Sklyanin Lie cobracket. \end{theorem}

\begin{proof} Let $\complex$ be spanned by $x^{-}$ and its dual $\dual{\complex}$ be spanned by $x^{+}$.  As braided-Lie bialgebras, $\complex$ and $\dual{\complex}$ are trivial: they have zero bracket and braided cobracket, by anti-symmetry.  Note that we can therefore dispense with the ``op'' on $\dual{\complex}$.  We fix the dual pairing as $\ip{x^{-}}{x^{+}}=1$.

We have made a central extension to $\Lie{g}$: explicitly, let this be $\widetilde{\Lie{g}}=\Lie{g} \dsum \complex\cdot h$.  The central extension acts on $\complex$ by $h \laction x^{-}=x^{-}$ and this induces $h \laction x^{+} = -x^{+}$ on the dual.  This centrally-extended algebra becomes a bialgebra with quasitriangular structure $\tilde{r}=r+h\tensor h$.

We now make the double-bosonisation and examine the resulting brackets.  Firstly, $\Lbracket{\Lie{g}}{x^{-}}=\Lbracket{\Lie{g}}{x^{+}}=0$ since $\Lie{g}$ is acting trivially and $h$ spans a central extension: $\Lbracket{\Lie{g}}{h}=0$.  To see that we have a copy of $\Lie{sl}_{2}(\complex)$ from $\complex\cdot x^{-} \dsum \complex\cdot h \dsum \complex\cdot x^{+}$, we must calculate the bracket between $x^{-}$ and $x^{+}$ as given by the double-bosonisation formul\ae.  We have \begin{align*} \Lbracket{x^{-}}{x^{+}} & = (x^{-})_{\underline{(1)}}\ip{x^{+}}{(x^{-})_{\underline{(2)}}} + (x^{+})_{\underline{(1)}}\ip{(x^{+})_{\underline{(2)}}}{x^{-}} \\ & \qquad + 2\tilde{r}^{(1)}_{+}\ip{x^{+}}{\tilde{r}^{(2)}_{+} \laction (x^{-})} \\ & = 0 + 0 + 2r^{(1)}_{+}\ip{x^{+}}{r^{(2)}_{+} \laction (x^{-})}+ 2h\ip{x^{+}}{h\laction x^{-}} \\ & = 0 + 2h\ip{x^{+}}{x^{-}} \\ & = 2h. \end{align*}  We can re-choose our basis vectors as $H=-2h$, $X^{-}=x^{-}$, $X^{+}=x^{+}$ to see that we indeed have $\Lie{sl}_{2}(\complex)$ as a Lie algebra.

So, we have $\dbos{\complex}{\widetilde{g}}{\op{\dual{\complex}}} \iso \Lie{g}\dsum \Lie{sl}_{2}(\complex)$ as Lie algebras and it remains to check that we have a direct sum as Lie coalgebras.  Double-bosonisation gives us the Lie cobracket on $\complex$ as follows. \begin{align*} \delta X^{-} = \delta x^{-} & = \underline{\delta} x^{-} + \tilde{r}^{(2)} \tensor \tilde{r}^{(1)} \laction x^{-} - \tilde{r}^{(1)} \laction x^{-} \tensor \tilde{r}^{(2)} \\ & = 0 + {r}^{(2)} \tensor {r}^{(1)} \laction x^{-} - {r}^{(1)} \laction x^{-} \tensor {r}^{(2)} \\ & \qquad + h\tensor h\laction x^{-} - h\laction x^{-} \tensor h \\ & = 0 + h\tensor x^{-} - x^{-}\tensor h \\ & = \frac{1}{2}(X^{-}\tensor H - H\tensor X^{-}) = \frac{1}{2}(X^{-} \wedge H). \end{align*}  Similarly, we have $\delta X^{+} = \frac{1}{2}( X^{+} \wedge H)$.  Equivalently, we see that the quasitriangular structure given by double-bosonisation (Proposition~\ref{doubleqtstr}) is \[ r^{\mathrm{new}}=\tilde{r}+\sum_{a} f^{a} \tensor e_{a} = r+h\tensor h+x^{+} \tensor x^{-}=r+\frac{1}{4} H\tensor H + X^{+} \tensor X^{-} \] where $\sum_{a} f^{a} \tensor e_{a}$ is a sum over $\{e_{a} \}$ a basis for $\Lie{b}=\complex$ and $\{ f^{a} \}$ is a dual basis.  This is the Drinfel\cprime d--Sklyanin quasitriangular structure.  Hence, we have a direct sum as bialgebras.
\end{proof}

\noindent We note that this result is independent of the choice of quasitriangular structure on $\Lie{g}$, since $\Lie{g}$ acts trivially in any case.

With respect to the induction procedure, this result excludes the choice $\Lie{b}_{-1}=V(\omega_{0})=\complex$ for all simple Lie algebras.  Note that choosing $\Lie{b}_{-1}=\complex$ fixes $\Lie{b}_{j}=0$ for all $j\leq -2$ by property~(\ref{suboftensorprodcond}) above together with the anti-symmetry required by a (graded) Lie bracket.  We have $\Wedge{2}{\complex}=0$ and $\Lie{b}_{-2}$ is required to be a submodule of this, so is zero, and this forces all the remaining $\Lie{b}_{j}$ to be zero.  So, to our list we add the property \begin{property} \Lie{b}_{-1}\ \mbox{is not trivial.} \end{property} {}\vspace{-1em}

\subsection{Classification of defining modules}\label{defmods}

We now classify the irreducible defining modules for the simple Lie algebras, that is, those irreducible modules satisfying conditions (\ref{1dwtspcond}) and (\ref{2domwtscond}) above.  To do this, we combine a result of Howe (as described by Stembridge (\cite{Stembridge})) with some analysis of dominant weights.  Howe's result classifies weight-multiplicity-free highest weight modules, that is, those with all weight spaces associated to non-zero weights being one-dimensional.  This is almost property (\ref{1dwtspcond}) above.  We then examine this relatively short list to determine the defining modules for each simple Lie algebra.

In our notation and the terminology of Stembridge, Howe's result is as follows:

\begin{theorem}[{\cite{Howe}}] Let $\Lie{g}$ be a finite-dimensional simple complex Lie algebra.  Then a non-trivial irreducible $\Lie{g}$-module $V(\lambda)$ has one-dimensional weight spaces if and only if \begin{enumerate} \item $\lambda$ is minuscule, \item $\lambda$ is quasi-minuscule and $\Lie{g}$ has only one short simple root, \item $\Lie{g}=C_{3}=\Lie{sp}_{6}$ and $\lambda=\omega_{3}$, or \item $\Lie{g}=A_{l}=\Lie{sl}_{l+1}$ and $\lambda=m\omega_{1}$ or $\lambda=m\omega_{l}$ for some $m\in \nat$. \end{enumerate} \end{theorem} \vspace{0.7em}

A weight $\lambda$ is called \emph{minuscule} if $\ip{\lambda}{\alpha}\leq 1$ for all $\alpha\in R$.  In \cite{Humphreys}, a dominant minuscule weight is called minimal and an alternative characterisation is given, namely, if $\mu$ is also dominant and $\mu \prec \lambda$ then $\mu=\lambda$.  Here $\prec$ is the usual partial ordering on weights.  In \cite{PSV}, non-zero minuscule dominant weights are called microweight.  We include the zero weight in the minuscule weights.  Minuscule weights are also discussed in the exercises to Chapter 6 of \cite{Bourbaki}.  Note that non-zero minuscule weights do not exist for all Dynkin types.

A weight $\lambda$ is called \emph{quasi-minuscule} if $\ip{\lambda}{\alpha}\leq 2$ for all $\alpha \in R$ and $\ip{\lambda}{\alpha'}=2$ for a unique $\alpha' \in R$.  For an irreducible root system, there is a unique dominant quasi-minuscule weight, namely the short dominant root.  The modules $V(\lambda)$, $\lambda$ quasi-minuscule, are called short-root representations in \cite{PSV}.  

Table 2 in \cite{PSV} gives the following lists of (non-zero) minuscule and quasi-minuscule weights: 
\begin{description}
\item[Non-zero minuscule weights:] \hfill

$\begin{array}{ll}  	
			A_{l} & \omega_{i},\ 1\leq i \leq l \\ 
			B_{l} & \omega_{l} \\ 
			C_{l} & \omega_{1} \\ 
			D_{l} & \omega_{1},\ \omega_{l-1},\ \omega_{l} \\ 
			E_{6} & \omega_{1},\ \omega_{6} \\ 
			E_{7} & \omega_{7} 
\end{array} $
\item[Quasi-minuscule weights:] \hfill

$ \begin{array}{ll} 
\begin{array}[t]{ll}  	A_{l} & [1,0,0,\ldots,0,1]\ (\mbox{adjoint}) \\ 
			B_{l} & \omega_{1} \\ 
			C_{l} & \omega_{2} \\ 
			D_{l} & \omega_{2} \\ 
			E_{6} & \omega_{2} 
\end{array} 
&
\begin{array}[t]{ll}	E_{7} & \omega_{1} \\ 
			E_{8} & \omega_{8} \\ 
			F_{4} & \omega_{4} \\ 
			G_{2} & \omega_{1} 
\end{array}
\end{array} $
\end{description}

The modules satisfying (\ref{1dwtspcond}) and (\ref{2domwtscond}) (the defining modules) can therefore be calculated by striking out of the above classification all those with too many orbits.

\begin{theorem}[{\cite{Bazlov}}]\label{defmodsthm} Let $\Lie{g}$ be a finite-dimensional simple complex Lie algebra.  The following is a list of all weights $\lambda$ such that the highest weight $\Lie{g}$-module $V(\lambda)$ satisfies properties (\ref{1dwtspcond}) and (\ref{2domwtscond}).

\vspace{1em}
$\begin{array}{ll}
\begin{array}[t]{ll} 
A_{1} & \omega_{0},\ \omega,\ 2\omega,\ 3\omega \\
A_{l},\ l\geq 2 & \omega_{0},\ \omega_{i}\ (1\leq i \leq l),\ 2\omega_{1},\ 2\omega_{l} \\
B_{l} & \omega_{0},\ \omega_{1},\ \omega_{l} \\
C_{3} & \omega_{0},\ \omega_{1},\ \omega_{3} \\
C_{l},\ l\geq 4 & \omega_{0},\ \omega_{1} \\
D_{l} & \omega_{0},\ \omega_{1},\ \omega_{l-1},\ \omega_{l}
\end{array} 
&
\begin{array}[t]{ll}
E_{6} & \omega_{0},\ \omega_{1},\ \omega_{6} \\
E_{7} & \omega_{0},\ \omega_{7} \\ 
E_{8} & \omega_{0} \\ 
F_{4} & \omega_{0} \\
G_{2} & \omega_{0},\ \omega_{1}
\end{array}
\end{array}$
\vspace{1em}
\end{theorem}
	
\begin{proof} The trivial module $V(\omega_{0})$ satisfies (\ref{1dwtspcond}) and (\ref{2domwtscond}) for all types.  It is well-known that the minuscule weights give rise to modules with exactly one Weyl group orbit (\cite{Humphreys}, \cite{PSV}).  Indeed, this is often given as essentially the definition.

Taking the types with only one short simple root excludes $\omega_{2}$ for $C_{l}$ of the quasi-minuscule weights since $C_{l}$ has $l-1$ short simple roots and to avoid repetitions in the labelling we have $l\geq 3$.  Of the remaining quasi-minuscule weights, we exclude the algebra-weight pairs corresponding to adjoint representations, namely $(A_{l}\ (l\geq 2), [1,0,0,\ldots,0,1])$, $(D_{l},\omega_{2})$, $(E_{6},\omega_{2})$, $(E_{7}, \omega_{1})$ and $(E_{8},\omega_{8})$, since in these cases the zero weight occurs with multiplicity $l$, the rank of $\Lie{g}$, which is greater than one.  

For $(F_{4},\omega_{4})$, the zero weight has multiplicity two, so is excluded.
We find that $V(\omega_{3})$ for $C_{3}$ (the long root representation) has two Weyl group orbits; the zero weight does not occur.  (We used LiE (\cite{LiE}) to obtain this information.)

For $(A_{1}, m\omega\ (=m\omega_{1}=m\omega_{l}))$, we have $\omega$ already---it is minuscule---and if $m\geq 4$, it is easy to see that $V(m\omega)$ has more than two orbits.  We are left with $m=2\ \mbox{or}\ 3$.  For $m=2$, we have $V(2\omega)=\Sym{2}{V}$ ($V$ the natural representation) and this has two orbits, the zero weight orbit and one other.  For $m=3$, $V(3\omega)=\Sym{3}{V}$ does not contain the zero weight but does have exactly two orbits.

Finally, for $(A_{l}\ (l\geq 2), m\omega_{1})$, $m=1$ is covered by the minuscule case and if $m\geq 3$ there are more than two orbits, as is easily seen.  However, $(A_{l}\ (l\geq 2), 2\omega_{1})$ is kept: $V(2\omega_{1})=\Sym{2}{V}$ ($V$ the natural representation) has exactly two orbits.  Since $m\omega_{l}$ is dual to $m\omega_{1}$, we are done. \end{proof}

\subsection{Induction for the exceptional series}\label{indandexcep}

In the remainder, we examine the question of extending the (known) exceptional series.  In particular, we show how our method indicates the obstructions to there being a simple $E_{9}$, $F_{5}$ or $G_{3}$.  We will see that two situations occur.  The first is that there may be no appropriate choices of modules to feed into the induction, as a result of the classification of the previous section.  The second is a lack of consistency, as described further below.

Our general algorithm is as follows, suggested by the six properties we listed in Section~\ref{necccond}.  We should take a simple algebra $\smun{\Lie{g}}{0}$ of rank $l$ and examine the list of defining modules in Theorem~\ref{defmodsthm} to find a candidate $V(\lambda_1)$ for $\Lie{b}_{-1}$, the first graded part of the braided-Lie bialgebra $\Lie{b}$ we need.  By Theorem~\ref{nottrivial}, we exclude the trivial module $V(\omega_{0})$ as a choice for $\Lie{b}_{-1}$.  Next calculate $\Wedge{2}{(V(\lambda_1))}$: if this is zero or has no irreducible submodules which are defining modules (for $\smun{\Lie{g}}{0}$), we stop here.  Otherwise, such a submodule, together with the zero subspace, is a candidate $V(\lambda_2)$ for $\Lie{b}_{-2}$.  We then see if there are non-zero maps from $\Lie{b}_{-1} \tensor \Lie{b}_{-2}$ into any defining module $V(\lambda_3)$ for $\smun{\Lie{g}}{0}$ satisfying the properties of a bracket, namely anti-symmetry and the (graded) Jacobi identity.  If there is such a map, then $V(\lambda_3)$ is a candidate for $\Lie{b}_{-3}$, and we repeat the process, considering maps from $\Lie{b}_{j} \tensor \Lie{b}_{k}$ to defining modules to find candidates for $\Lie{b}_{j+k}$.

We now apply this algorithm to the appropriate simple algebras of rank 8, 4 and 2, to attempt to construct $E_{9}$, $F_{5}$ and $G_{3}$.

\subsubsection{$E_{9}$}

The first obvious line of attack is to consider induction from $E_{8}$.  We may deal with this easily, as Theorems \ref{nottrivial} and \ref{defmodsthm} show that there are no possible choices for $\Lie{b}_{-1}$ and hence no inductions.  In fact, this is a stronger statement than we need as by considering the deletion from $E_{9}$, we would require $\Lie{b}_{-1}=V(\omega_{8}; E_{8})=E_{8}$ (the adjoint representation).  Clearly, we cannot have this as there is the eight-dimensional Cartan subalgebra, so the zero weight space is not one-dimensional.

Of course, we could instead look to induce from the other series.  If we consider induction from $D_{8}$, then we will require $\Lie{b}_{-1}=V(\omega_{7}; D_{8})$; the embedding of $D_{8}$ in $E_{9}$ we choose is $\iota:i\mapsto 10-i$, where the labelling of the Dynkin diagram for $E_{9}$ follows the usual pattern for $E_{l}$, $l=6,7,8$.  As desired, $V(\omega_{7}; D_{8})$ is a defining module for $D_{8}$ but we have $\Wedge{2}{(V(\omega_{7}; D_{8}))}=V(\omega_{2}; D_{8}) \dsum V(\omega_{6}; D_{8})$ and neither of these is defining.  So we are forced to take $\Lie{b}_{j}=0$ for $j \leq -2$.  Thus our candidate for $E_{9}$ is $\Lie{g}=V(\omega_{7}; D_{8}) \rtimesdot \widetilde{D_{8}} \ltimesdot \dual{V(\omega_{7}; D_{8})}$, which has dimension 377.  In this case, we have only the zero braided-Lie bialgebra structure and note that Corollary~4.2 of \cite{BraidedLie} does not apply here (to tell us $\Lie{g}$ is simple) since $\Wedge{2}{V(\omega_{7}; D_{8})}$ is not isotypical.

From $A_{8}$, the situation is more complicated.  To give the correct Dynkin diagram, we must choose $\Lie{b}_{-1}=V(\omega_{3}; A_{8})$; the embedding is \[ \iota= \left( \begin{array}{llllllll} \Ss 1 & \Ss 2 & \Ss 3 & \Ss 4 & \Ss 5 & \Ss 6 & \Ss 7 & \Ss 8 \\ \Ss 1 & \Ss 3 & \Ss 4 & \Ss 5 & \Ss 6 & \Ss 7 & \Ss 8 & \Ss 9 \end{array} \right). \]  Now $\Wedge{2}{V(\omega_{3}; A_{8})}=V([0,1,0,1,0,0,0,0]; A_{8}) \dsum V(\omega_{6}; A_{8})$ and $V(\omega_{6}; A_{8})$ is a defining module for $A_{8}$ so we have the choice of $V(\omega_{6}; A_{8})$ and the zero space for $\Lie{b}_{-2}$.  Next, we have (all as $A_{8}$-modules) \begin{multline*} V(\omega_{3}) \tensor V(\omega_{6})= V([0,0,1,0,0,1,0,0]) \dsum V([0,1,0,0,0,0,1,0]) \\ \dsum V([1,0,0,0,0,0,0,1]) \dsum V(\omega_{0}) \end{multline*} and the first three terms are not defining but the last is.  (We excluded $V(\omega_{0})$ as a choice for $\Lie{b}_{-1}$ in Theorem~\ref{nottrivial} but it is valid as a choice for other $\Lie{b}_{j}$ and indeed it does occur.)  Hence we have the choices $\Lie{b}_{-3}=V(\omega_{0}; A_{8})$ or $\Lie{b}_{-3}=0$.  For $\Lie{b}_{-4}$, we have $V(\omega_{3}) \tensor V(\omega_{0})=V(\omega_{3})$ and \[ \Wedge{2}{(V(\omega_{6}))}=V([0,0,0,0,1,0,1,0]) \dsum V(\omega_{3}) \] so we can choose either $V(\omega_{3}; A_{8})$ or zero for $\Lie{b}_{-4}$.  Finally, for $\Lie{b}_{-5}$ and higher parts, we will see the same pattern, namely \[ \Lie{b}_{j} = \begin{cases} V(\omega_{3}; A_{8}) & \text{if}\ j\equiv -1 \!\!\mod{3} \\ V(\omega_{6}; A_{8}) & \text{if}\ j\equiv -2 \!\!\mod{3} \\ V(\omega_{0}; A_{8}) & \text{if}\ j\equiv 0 \!\!\mod{3}. \end{cases} \]

Observe that \begin{align*} & \dim \left( V(\omega_{3}; A_{8}) \rtimesdot \widetilde{A_{8}} \ltimesdot \dual{V(\omega_{3}; A_{8})} \right) = 249\quad \text{and} \\ &  \dim \left( \left(V(\omega_{3}; A_{8}) \dsum V(\omega_{6}; A_{8})\right) \rtimesdot \widetilde{A_{8}} \ltimesdot \dual{\left(V(\omega_{3}; A_{8}) \dsum V(\omega_{6}; A_{8})\right)} \right) = 417. \end{align*}  Recall that the proposed candidate for $E_{9}$ found above by induction from $D_{8}$ had dimension 377.  These are clearly inconsistent, so even discounting the lack of an induction from $E_{8}$, we have no sensible candidate for $E_{9}$ which agrees from both the $A$ and $D$ inductions.  

This analysis suggests that it is unlikely that in any larger scheme of finite-dimensional algebras we would find even a semisimple $E_{9}$ candidate.  Of course, we know of an (infinite-dimensional) $E_9$: the affine Kac-Moody algebra.  In the above scenario we would have to accept that no $E_{9}$ can exist as soon as we know that we cannot reach it from $E_{8}$.  Then this instantly rules out any diagram with the diagram $E_{9}$ as a sub-diagram, of course.  In a similar way, this reduces the number of cases to be considered in a proof of the classification of the simples significantly.

\subsubsection{$F_{5}$}

As for $E_{9}$, Theorems \ref{nottrivial} and \ref{defmodsthm} exclude the possibility of an induction from the natural starting point $F_{4}$, giving rise to the candidate for $F_{5}$ with Dynkin diagram \[ \begin{array}{c}\scalebox{0.4}{\includegraphics{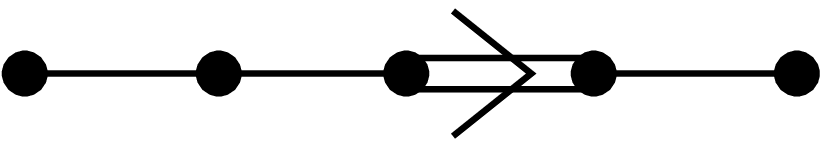}} \end{array} \]  For the other choice, namely $\Lie{b}_{-1}=V(\omega_{4}; F_{4})$, we simply use Theorem~\ref{defmodsthm}.

From the remaining series, we may start from either $B_{4}$ or $C_{4}$.  Similar arguments to the above rule out any induction from $C_{4}$ but from $B_{4}$, we obtain the candidate $V(\omega_{4}; B_{4}) \rtimesdot \widetilde{B_{4}} \ltimesdot \dual{V(\omega_{4}; B_{4})}$ of dimension 69 corresponding to the diagram given above.  In this case, we only have a single candidate but since $\dim F_{4}=52$, for consistency we would need an $F_{4}$-module of dimension $(69-52-1)\div 2=8$.  Such a module does not exist.

\subsubsection{$G_{3}$}

For $G_{3}$, we may consider adding the new node to either the first node in $G_{2}$ or the second.  By Theorem~\ref{defmodsthm}, we cannot add it to the first, as this would require $\Lie{b}_{-1}=V(\omega_{2}; G_{2})$ and this is not a defining module for $G_{2}$.  However, $V(\omega_{1}; G_{2})$ is a defining module so we may choose $\Lie{b}_{-1}=V(\omega_{1}; G_{2})$.  Then $\Wedge{2}{(V(\omega_{1}; G_{2}))}=V(\omega_{1}; G_{2}) \dsum V(\omega_{2}; G_{2})$, so we may choose $\Lie{b}_{-2}=0$ or $\Lie{b}_{-2}=V(\omega_{1}; G_{2})$.  If we choose the latter, we have an appropriate map to allow us to choose $\Lie{b}_{-3}=V(\omega_{1}; G_{2})$ and so on.  Clearly, we cannot go on choosing $V(\omega_{1}; G_{2})$ forever so we must decide if \begin{align}\label{g2g3} \left( \bigdsum_{j=-1}^{-m} V(\omega_{1}; G_{2}) \right) \rtimesdot \widetilde{G_{2}} \ltimesdot \op{ \dual{ \left( \bigdsum_{j=-1}^{-m} V(\omega_{1}; G_{2}) \right)}} \end{align} is simple for some value of $m$, with the appropriate braided-Lie bialgebra structure on the direct sum.  Now, we know that this double-bosonisation cannot be simple but it is not immediately obvious why not.  This case again illustrates that the list of properties in Section~\ref{necccond} is not yet complete.

We may only obtain $G_{3}$ from $A_{2}$, other than from $G_{2}$, but we may do this in two different ways (compare with the two embeddings of $A_{1}$ in $G_{2}$), leading to the possible diagrams \[ \begin{array}{ccc} \scalebox{0.4}{\includegraphics{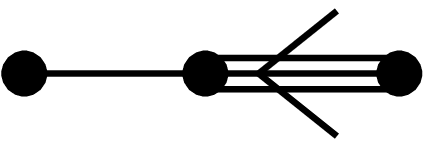}} \qquad & \text{and} & \qquad \scalebox{0.4}{\includegraphics{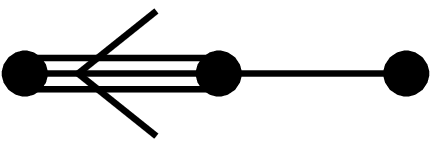}} \end{array} \]  We may exclude the first of these by the usual appeal to Theorem~\ref{defmodsthm}: $V(3\omega_{i}; A_{2})$ is not defining for either $i=1$ or $i=2$.  For the second, we find ourselves in a similar periodic situation to that for $E_{9}$, with  \[ \Lie{b}_{j} = \begin{cases} V(\omega_{1}; A_{2}) & \text{if}\ j\equiv -1 \!\!\mod{3} \\ V(\omega_{2}; A_{2}) & \text{if}\ j\equiv -2 \!\!\mod{3} \\ V(\omega_{0}; A_{2}) & \text{if}\ j\equiv 0 \!\!\mod{3}. \end{cases} \]  Some calculations with dimensions show that we do not have the same consistency problem between this induction and that discussed above from $G_{2}$.  For example, choosing $m=2$ in (\ref{g2g3}) gives an algebra of dimension $43$, which is consistent with choosing $\Lie{b}_{-j}$ subject to the above rule and non-zero for $j=1,\ldots, 7$.  Furthermore, since \[ \dim \left( V(\omega_{1}; A_{2}) \dsum V(\omega_{2}; A_{2}) \dsum V(\omega_{0}; A_{2}) \right) = 7, \] we can find similar matching candidates for each choice of $m$.  This appears to be as far as we can progress with this case by considering just the module structures.

\section{Concluding remarks}\label{concl}

As we have described in the introduction, this work is the start of a programme which aims to use Lie induction to provide insight into the simple Lie (bi-)algebras.  There are some outstanding questions prompted even by the small number of examples we have given, which we record here.

We have described an algorithm for working out inductions at the start of Section~\ref{indandexcep}.  However, we have seen in that section that this process may not terminate.  Furthermore, we have an example where although it does indeed terminate, the resulting braided-Lie bialgebra does not give rise to a simple Lie algebra.  In short, we can be sure that additional properties to the six given in Section~\ref{necccond} are necessary to produce a classification of the simple algebras.  There appear to be a number of questions still to be answered before such a proof could be produced, most notably the following: \begin{enumerate} 
\item How do we fix the remaining values in the Cartan matrix for the induced algebra $\Lie{g}$?  Note that we have the sub-matrix corresponding to $\smun{\Lie{g}}{0}$ and the row corresponding to the new simple root from the highest weight of $\Lie{b}_{-1}$.  
\item Where, algebraically, do the properties of the Cartan matrix come from, for example, $C_{ij}=0 \iff C_{ji}=0$?  The restriction on the values in the Cartan matrix is related to the restriction on the number of non-zero graded parts we may have in $\Lie{b}$---where do these come from?
\item Is there a general and/or easy test to decide if a double-bosonisation 
$\dbos{b}{\widetilde{\smun{g}{\mathrm{0}}}}{\op{\dual{b}}}$ is simple?
\end{enumerate}

\noindent These are clearly not independent: an answer to the first two questions would give us an effective answer to the third.  There are also some wider questions of consistency, some of which we have touched on during our worked examples.  For example, for any Dynkin diagram corresponding to a simple Lie algebra (in particular, this is connected), any connected sub-diagram corresponds to a simple Lie subalgebra.  Clearly, we need this property to hold with respect to induction, or else our algebras are not well-defined.

We have here built on the work begun in \cite{BraidedLie} and many of the comments concluding that paper apply equally well here.  We have also restricted ourselves to working over the complex field and to considering the standard quasitriangular structure.  As noted there, it ought to be possible consider twisting and $*$-structures on the braided-Lie bialgebras and we now see that they would have to be compatible with the graded structure.

We also recall that the double-bosonisation construction can be defined when working over fields of any characteristic except two.  The theorem of Azad, Barry and Seitz (\cite{ABS}) we use in Section~\ref{grading} holds except for the following algebra-characteristic pairs, called \emph{special}: $(B_{l},\chr K=2)$, $(C_{l},2)$, $(F_{4},2)$, $(G_{2},2)$ and $(G_{2},3)$.  So we expect that the inductive method ought to carry over to (most) positive characteristics.  This, and an analysis of these special pairs, would be an interesting direction for further work.

A major motivation for this work was the hope that analysing this construction for Lie bialgebras would provide insight into the corresponding (known) construction for Hopf algebras, particularly the quantum groups $\qea{g}$.  Na\"{\i}vely, we expect a close relationship, especially given the similarity of the representation theories.  We hope to develop this elsewhere.

Finally, we have dealt here with the finite-dimensional case only but it seems sensible to extend our field of view to Kac-Moody Lie algebras in general.  The definitions of a braided-Lie bialgebra and of double-bosonisation do not require finite-dimensionality: the only result we use that does is the quasitriangularity of the double-bosonisation but with care this should not be a problem.  We may require formal power series, for example.  We intend to consider induction from the finite-type simple algebras to the affine ones in detail elsewhere, also.

%Built-in appendix command doesn't work!?

\renewcommand{\thesection}{\Alph{section}}
\setcounter{section}{0}
\section{Appendix}\label{appendix}

\subsection{Calculations: some tools}

We now proceed to the explicit calculations for the deletions $(\Lie{g}, d, \smun{\Lie{g}}{0}, \iota)$ with $\Lie{g}$ and $\smun{\Lie{g}}{0}$ simple.  We first set our notations and indicate how we may use the results of Section~\ref{deletion} to simplify the calculations.

We will give $\iota$ in one of three forms.  Recall that $\iota$ is an embedding of $\smun{\Lie{g}}{0}$ into $\Lie{g}$ but that this is equivalent to an embedding of Dynkin diagrams, and hence can be expressed as a map of the labels of the diagram nodes.  We write $\id$ if $\iota$ is the identity map or write $\iota$ algebraically, if possible.  For example, we may write $i \mapsto i+1$ for the embedding of (the diagram) $A_{3}$ into $A_{4}$ given by $1\mapsto 2$, $2\mapsto 3$, $3\mapsto 4$.  Otherwise we will write $\iota$ in two-row permutation notation, with $\smun{\Lie{g}}{0}$ on top, although it will not be a genuine permutation as the label sets will differ.

In the notation of Section~\ref{grading}, set $\Lie{b}_{i}=\smun{\Lie{g}}{[i]}$ for $i<0$, the graded components of $\Lie{b}$ as a graded Lie algebra.  The grading gives us another way to analyse $\Lie{b}$, since we can consider the $\smun{\Lie{g}}{0}$-module $\Wedge{2}{\Lie{b}}$ and its subspaces.  In particular, we can consider $\Wedge{2}{\Lie{b}_{-1}}$, which will give us information about $\Lie{b}_{-2}$.  \label{submods}

Firstly, if $m_{d}=1$, so $\Lie{b}$ is irreducible, $\Lie{b}$ has zero bracket.  Secondly, if $m_{d}=2$ and $\dim \Lie{b}_{-2}=1$, there is a non-zero bracket on $\Lie{b}_{-1}$ and it is a cocycle central extension of the zero bracket.  For $\dim \Lie{b}_{-2}=1$ implies $\Lie{b}_{-2}$ is spanned by $\Lambda$, the highest root in $\Lie{g}$ and the grading on $\Lie{b}$ tells us that if $X_{\alpha}^{-},\ X_{\beta}^{-} \in \Lie{b}_{-1}$ then $\Lbracket{X_{\alpha}^{-}}{X_{\beta}^{-}}= \delta(\alpha + \beta, \Lambda) X_{\Lambda}^{-}$ where $\delta(\alpha + \beta,\Lambda)=0$ if $\alpha+\beta \neq \Lambda$ and $\delta(\alpha+\beta,\Lambda)=c_{\alpha \beta}$ (some constant depending on $\alpha$ and $\beta$) if $\alpha+\beta = \Lambda$.  If $m_{d} \geq 2$ and $\dim \Lie{b}_{-2} > 1$, although a similar additive formula will hold, we cannot be so explicit.

The bracket $\Lbracket{\ }{\ }_{\Lie{g}} : \Wedge{2}{\Lie{g}} \to \Lie{g}$ clearly restricts to $\Lbracket{\ }{\ }_{\Lie{b}} : \Wedge{2}{\Lie{b}} \to \Lie{b}$ and indeed even restricts to $\Lbracket{\ }{\ }_{-1} : \Wedge{2}{\Lie{b}_{-1}} \to \Lie{b}_{-2}$.  Hence we can consider the kernel $K_{-1}$ of $\Lbracket{\ }{\ }_{-1}$, which must be a sum of irreducible components of $\Wedge{2}{\Lie{b}_{-1}}$ (possibly zero but not all of $\Lie{b}_{-1}$) and so we have $\Lie{b}_{-2} \iso \Wedge{2}{\Lie{b}_{-1}} / K_{-1}$.  Given the restricted number of possibilities for $\Lie{b}_{-1}$ (which we know), clearly there will not be very many choices for $\Lie{b}_{-2}$, so in the case where $\Lie{g}$ and $\smun{\Lie{g}}{0}$ are simple ($m_d \leq 3$) we are essentially done.  In particular, if $\Lie{b}_{-2} \neq \zeroset$ and $\Wedge{2}{\Lie{b}_{-1}}$ is irreducible, we have $\Ker{\Lbracket{\ }{\ }_{-1}} \neq \Wedge{2}{\Lie{b}_{-1}}$ so $\Ker{\Lbracket{\ }{\ }_{-1}}=0$ and $\Lie{b}_{-2} \iso \Wedge{2}{\Lie{b}_{-1}}$.

All of the above has been classical, in the sense that it has been derived from properties of root systems.

We now consider the final structure we need on $\Lie{b}$, that of a braided-Lie bialgebra.  This has been given in the proof of \cite[Proposition 4.5]{BraidedLie} when the quasitriangular structure on $\Lie{g}$ is chosen to be the Drinfel\cprime d--Sklyanin solution.  In this case, it has the general form \[ \underline{\delta} X_{\alpha}^{-} = \sum_{\alpha = \beta + \gamma} c_{\beta \gamma} X_{\beta}^{-} \wedge X_{\gamma}^{-} \in \Wedge{2}{\Lie{b}}. \]

By the additivity property of the multiplicity $\mult{d}{-}$, this must be zero on elements of $\Lie{b}_{-1}$ since if $\alpha = \beta+\gamma$ for some $\beta,\ \gamma \in \Lie{b}_{-1}$ then $\mult{d}{\alpha} = 2$.  However, if $m_{d} \geq 2$, $\underline{\delta}$ will not be zero on $\Lie{b}_{j}$, $j \leq -2$.  If $\Wedge{2}{\Lie{b}_{-1}}$ is irreducible, by the above, $\Lie{b}_{-2} \iso \Wedge{2}{\Lie{b}_{-1}}$ so using Schur's lemma $\underline{\delta}$ must be an isomorphism.

We have used the above tools and the computer program LiE (\cite{LiE}) to calculate the braided-Lie bialgebra structures arising in the deletions $(\Lie{g},d,\smun{\Lie{g}}{0},\iota)$ for all choices of $\Lie{g}$ and $d$ such that $\Lie{g}$ and $\smun{\Lie{g}}{0}$ are simple.  These calculations are given below, grouped by the value $m_{d}$ for each deletion.  For the exceptional simple Lie algebras, we have given less detail as the maps are not easily expressible in simple terms and the explicit formul\ae\/ not necessarily very informative.  We wish to stress, though, that once the task of writing down the Weyl basis (or equivalently the root system) has been achieved, it is relatively simple to recover these formul\ae.  For a summary of the module structures, we refer the reader to Table~\ref{deletionstable} on page~\pageref{deletionstable}.

\subsection{$m_{d}=1$}\label{md1}

Recall from above that in the case $m_{d}=1$, $\Lie{b}=\Lie{b}_{-1}$ is irreducible and has zero Lie algebra and braided-Lie coalgebra structures.  Below we give the induced isomorphisms of $\Lie{b}$ as a set of roots of $\Lie{g}$ with the usual basis for $\Lie{b}$ as a $\smun{\Lie{g}}{0}$-module of the appropriate highest weight.  Rather than numbering the cases, we will use a two-letter code corresponding to the Dynkin types of $\Lie{g}$ and $\smun{\Lie{g}}{0}$ (in that order), suppressing the rank as subscript where this is appropriate.

\begin{description}
\item[($AA$)] Deletion $(A_{l+1}, l, A_{l}, \id )$

$\Lie{b}$ has highest weight $\omega_{l}$ so is the natural representation of $A_{l}=\Lie{sl}_{l+1}$ on the vector space $V$ of dimension $l+1$.  A basis for $V$ is $\{ \smun{e}{i} \mid 1 \leq i \leq l+1 \}$ and the highest weight vector is $\smun{e}{1}$.  The corresponding $\smun{\Lie{g}}{0}$-module isomorphism is $\smun{e}{i} \mapsto \smun{X}{(l-i+2)\cdots (l)(l+1)}^{-}$.

\item[($BB$)] Deletion $(B_{l+1},1 , B_{l}, i\mapsto i+1 )$

$\Lie{b}$ has highest weight $\omega_{1}$ so is the natural representation of $B_{l}=\Lie{so}_{2l+1}$ on the vector space $V$ of dimension $2l+1$.  A basis for $V$ is given by $\{ \smun{e}{i} \mid 1 \leq i \leq 2l+1 \}$ and the highest weight vector is $\smun{e}{1}$.  The corresponding $\smun{\Lie{g}}{0}$-module isomorphism is \[ \begin{cases} \smun{e}{i} \mapsto \smun{X}{12\cdots i}^{-} & \text{for $1\leq i \leq l$} \\ (-1)^{l+i}\smun{e}{l+i-1} \mapsto \smun{X}{12\cdots (i-1)(i)(i)\cdots (l+1)(l+1)}^{-} & \text{for $2 \leq i \leq l+1$} \\ -\smun{e}{2l+1} \mapsto \smun{X}{12\cdots (l+1)}^{-}. \end{cases} \]

\item[($DD$)] Deletion $(D_{l+1},1 , D_{l}, i \mapsto i+1)$

$\Lie{b}$ has highest weight $\omega_1$ so is the natural representation of $D_{l}=\Lie{so}_{2l}$ on the vector space $V$ of dimension $2l$.  A basis for $V$ is $\{ \smun{e}{i} \mid 1 \leq i \leq 2l \}$ and the highest weight vector is $\smun{e}{1}$.  The corresponding $\smun{\Lie{g}}{0}$-module isomorphism is \[ \begin{cases} \smun{e}{i} \mapsto \smun{X}{12\cdots i}^{-} & \text{for $1 \leq i \leq l$} \\ (-1)^{l+i}\smun{e}{l+i-1} \mapsto \smun{X}{12\cdots (i-1)(i)(i)\cdots (l-1)(l-1)(l)(l+1)}^{-} & \text{for $2 \leq i \leq l-1$} \\ \smun{e}{2l-1} \mapsto \smun{X}{12\cdots (l+1)}^{-} \\ -\smun{e}{2l} \mapsto \smun{X}{12\cdots (l-1)(l+1)}^{-}. \end{cases} \]

\item[($E_{7}E_{6}$)] Deletion $(E_{7},7 , E_{6}, \id)$

$\Lie{b}$ has highest weight $\omega_6$ and is one of the dual pair of representations of $E_6$ of dimension $27$.  As discussed in \cite{Schafer} and \cite{Baez}, these come from the action of the group $E_6$, as a group of determinant-preserving linear transformations, on the exceptional Jordan algebra $\smun{\Lie{h}}{3}(\octo )$.

\item[($CA$)] Deletion $(C_{l+1}, l+1, A_{l}, i \mapsto l-i+1 )$

$\Lie{b}$ has highest weight $2\omega_1$ so is the symmetric square $\Sym{2}{V}$ with $V$ the $(l+1)$-dimensional natural representation of $A_{l}$.  A basis for $\Sym{2}{V}$ is $\{ \smun{e}{i}\smun{e}{j} \mid 1 \leq i \leq j \leq l+1 \}$, so the dimension of $\Sym{2}{V}$ is $\frac{1}{2}(l+1)(l+2)$, and the highest weight vector is $\smun{e}{1}^{2}$.  The corresponding $\smun{\Lie{g}}{0}$-module isomorphism is \[ \begin{cases} \smun{e}{i}^{2} \mapsto \smun{X}{(l-i+2)(l-i+2)\cdots (l)(l)(l+1)}^{-} & \text{for $1 \leq i \leq l+1$} \\ \smun{e}{i}\smun{e}{j} \mapsto \smun{X}{(l-j+2)(l-j+3)\cdots (l-i+1)(l-i+2)(l-i+2)\cdots (l)(l)(l+1)}^{-} & \text{for $i<j$}. \end{cases} \]

\item[($DA$)] Deletion $(D_{l+1}, l+1, A_{l}, i \mapsto l-i+1)$

$\Lie{b}$ has highest weight $\omega_2$ so is the second exterior power $\Wedge{2}{(V)}$ with $V$ the $l+1$-dimensional natural representation of $A_{l}$.  The dimension of $\Wedge{2}{(V)}$ is $\frac{1}{2}l(l+1)$.  A basis for $\Wedge{2}{(V)}$ is $\{ \smun{e}{i} \wedge \smun{e}{j} \mid 1 \leq i < j \leq l+1 \}$ and the highest weight vector is $\smun{e}{1} \wedge \smun{e}{2}$.  The corresponding $\smun{\Lie{g}}{0}$-module isomorphism is \[ \left\{ \begin{array}{ll} \smun{e}{1} \wedge \smun{e}{2} \mapsto \smun{X}{l+1}^{-} & \\ \smun{e}{1} \wedge \smun{e}{j} \mapsto \smun{X}{(l-j+2)\cdots (l-1)(l+1)}^{-} & \qquad \qquad \quad \text{for $j \geq 3$} \\ \multicolumn{2}{l}{\smun{e}{i} \wedge \smun{e}{j} \mapsto \smun{X}{(l-j+2)\cdots(l-i+1)(l-i+2)(l-i+2)\cdots (l-1)(l-1)(l)(l+1)}^{-}} \\ & \qquad \qquad \quad \text{for $2\leq i < k \leq l+1$}. \end{array} \right. \]

\item[($E_{6}D_{5}$)] Deletion $(E_{6}, 1, D_{5}, i \mapsto 7-i)$

$\Lie{b}$ has highest weight $\omega_{4}$ so is the positive (half-)spin representation $S_{5}^{+}$ of $D_{5}$ (see for example \cite[Chapter 20]{FultonHarris}).  As a vector space, $S_{5}^{+}=\Wedge{0}{(V)}\dsum \Wedge{2}{(V)}\dsum \Wedge{4}{(V)}$ with $V$ the vector space of dimension 5.  Hence a basis for $S_{5}^{+}$ is given by taking the natural bases for these pieces.  The highest weight vector is $\smun{e}{1} \wedge \smun{e}{2} \wedge \smun{e}{3} \wedge \smun{e}{4}$.  The corresponding $\smun{\Lie{g}}{0}$-module isomorphism may easily be calculated from this.

\end{description}

\subsection{$m_{d}=2$}\label{md2}

The Lie algebra and braided-Lie coalgebra structures are no longer zero and we give explicit expressions for these where possible, in addition to the description following the pattern of the above.

\begin{description}
\item[($CC$)] Deletion $(C_{l+1}, 1, C_{l}, i\mapsto i+1 )$

$\Lie{b}_{-1}$ has highest weight $\omega_1$ so is the natural representation of $C_{l}=\Lie{sp}_{2l}$ on the vector space $V$ of dimension $2l$.  A basis for $V$ is $\{ \smun{e}{i} \mid 1 \leq i \leq 2l \}$ and the highest weight vector is $\smun{e}{1}$.  The corresponding $\smun{\Lie{g}}{0}$-module isomorphism is \[ \begin{cases} \smun{e}{i} \mapsto \smun{X}{12\cdots i}^{-} & \text{for $1\leq i\leq l$} \\ (-1)^{l+i-1}\smun{e}{l+i-1} \mapsto \smun{X}{12\cdots (i-1)(i)(i)\cdots (l)(l)(l+1)}^{-} & \text{for $2 \leq i \leq l$} \\ \smun{e}{2l} \mapsto \smun{X}{12\cdots (l+1)}^{-}. & \end{cases} \]

$\Lie{b}_{-2}$ has highest weight $\omega_0=[0,0,\ldots ,0]$ so is the trivial representation.  We can see this by a dimension calculation.  So, as described above, $\Lie{b}_{-2}$ is spanned by the highest root, $\smun{X}{1122\cdots (l)(l)(l+1)}^{-}=\varsigma$.

The bracket on $\Lie{b}=\Lie{b}_{-1} \dsum \Lie{b}_{-2}$ is a cocycle central extension of the zero bracket on $\Lie{b}_{-1}$, with $\Lbracket{\smun{e}{i}}{(-\smun{e}{l+i})}= c_{i}\,\varsigma$ for $1 \leq i \leq l-1$ and $\Lbracket{\smun{e}{l}}{\smun{e}{2l}}=c_{l}\,\varsigma$, where the $c_i$, $1\leq i \leq l$, are constants.  The braided-Lie cobracket is zero on elements of $\Lie{b}_{-1}$, as discussed previously, and \[ \underline{\delta}\varsigma = \sum_{i=1}^{l} \gamma_{i}\,( \smun{e}{i} \wedge \smun{e}{l+i}) \] for some constants $\gamma_i$.

We have $\Wedge{2}{\Lie{b}_{-1}} \iso V(\omega_{2}) \dsum V(\omega_{0})$ ($V(\omega)$ is the representation of $C_{l}$ with highest weight $\omega$) and we see that we have $\Ker{\Lbracket{\ }{\ }} \iso V(\omega_{2})$, $\Lie{b}_{-2} \iso V(\omega_{0}) = \complex$.

\item[($E_{8}E_{7}$)] Deletion $(E_{8}, 8, E_{7}, \id )$

$\Lie{b}_{-1}$ has highest weight $\omega_7$ and is the smallest non-trivial representation of $E_{7}$.  This may be realised by a Freudenthal triple system (see \cite{Baez} and the references therein).  The dimension of $\Lie{b}_{-1}$ is 56. 

$\Lie{b}_{-2}$ has highest weight $\omega_0$ so is the trivial representation, by a dimension calculation.  It is spanned by the highest root in $E_{8}$, $\smun{X}{(2,3,4,6,5,4,3,2)}^{-}$.

The bracket on $\Lie{b}=\Lie{b}_{-1} \dsum \Lie{b}_{-2}$ is again a cocycle central extension of the zero bracket on $\Lie{b}_{-1}$ and has the additive form described previously.  Similarly, the braided-Lie cobracket is non-zero only on $\Lie{b}_{-2}$ and has the additive form.

\end{description}

\begin{jegnote}{} One might consider that this deletion provides the most natural basis for the 56-dimensional representation of $E_{7}$. \end{jegnote}

\begin{description}
\item[($E_{6}A_{5}$)] Deletion $(E_{6}, 2, A_{5}, \left( \begin{array}{lllll} \Ss 1 & \Ss 2 & \Ss 3 & \Ss 4 & \Ss 5 \\ \Ss 1 & \Ss 3 & \Ss 4 & \Ss 5 & \Ss 6 \end{array} \right) )$

$\Lie{b}_{-1}$ has highest weight $\omega_3$ so is the third exterior power $\Wedge{3}{(V)}$ with $V$ the 6-dimensional natural representation of $A_{5}$.  The dimension of $\Wedge{3}{(V)}$ is 20.  A basis for $\Wedge{3}{(V)}$ is $\{ \smun{e}{i} \wedge \smun{e}{j} \wedge \smun{e}{k} \mid 1 \leq i < j < k \leq 6 \}$ and the highest weight vector is $\smun{e}{1} \wedge \smun{e}{2} \wedge \smun{e}{3}$.  The corresponding $\smun{\Lie{g}}{0}$-module isomorphism may be calculated from this.

$\Lie{b}_{-2}$ has highest weight $\omega_{0}$, so is the trivial representation, by a dimension calculation.  It is spanned by the highest root in $E_{6}$, $\smun{X}{(1,2,2,3,2,1)}^{-}$.  However, as we will see, we should consider $\Lie{b}_{-2}$ to be $\Wedge{6}{(V)}$ with $V$ as before, spanned by $\smun{e}{1} \wedge \smun{e}{2} \wedge \smun{e}{3} \wedge \smun{e}{4} \wedge \smun{e}{5} \wedge \smun{e}{6}$.

The bracket on $\Lie{b}=\Lie{b}_{-1} \dsum \Lie{b}_{-2}$ is given by the map $\wedge : \Lie{b}_{-1} \tensor \Lie{b}_{-1} \to \Lie{b}_{-2}$, $(\smun{e}{i_1} \wedge \smun{e}{j_1} \wedge \smun{e}{k_1}) \tensor (\smun{e}{i_2} \wedge \smun{e}{j_2} \wedge \smun{e}{k_2}) \mapsto \smun{e}{i_1} \wedge \smun{e}{j_1} \wedge \smun{e}{k_1} \wedge \smun{e}{i_2} \wedge \smun{e}{j_2} \wedge \smun{e}{k_2}$, that is, the wedge product.  The bracket is zero on all other elements of $\Lie{b}^{\tensor 2}$.  The braided-Lie cobracket is a map \[ \underline{\delta}: \Lie{b}_{-2} \to \Lie{b}_{-1} \wedge \Lie{b}_{-1} \iso \Wedge{3}{(V)} \wedge \Wedge{3}{(V)} \iso \Wedge{6}{(V)} \iso \Lie{b}_{-2} \] so must be a non-zero scalar multiple of the identity.

\item[($F_{4}C_{3}$)] Deletion $(F_{4},1,C_{3}, i\mapsto 5-i)$

$\Lie{b}_{-1}$ has highest weight $\omega_3$ and is described as the kernel of the contraction map $\phi_{3}:\Wedge{3}{(V)} \to V$ for $V$ the $6$-dimensional natural representation of $C_{3}=\Lie{sp}_{3}$ (see, for example, \cite[p.\:258]{FultonHarris}).  The dimension of $\Lie{b}_{-1}$ is 14.

$\Lie{b}_{-2}$ has highest weight $\omega_{0}$ so is the trivial representation, by a dimension calculation. It is spanned by the highest root in $F_{4}$, $\smun{X}{(2,3,4,2)}^{-}=\varsigma$.

The bracket on $\Lie{b}=\Lie{b}_{-1} \dsum \Lie{b}_{-2}$ is again a cocycle central extension of the zero bracket on $\Lie{b}_{-1}$ and has the additive form described previously.  Similarly, the braided-Lie cobracket is non-zero only on $\varsigma$ and has the additive form.

\item[($G_{2}A_{1}$)(a)] Deletion $(G_{2}, 2,A_{1}, \id )$

$\Lie{b}_{-1}$ has highest weight $3\omega_1$ so is the third symmetric power $\Sym{3}{V}$ with $V$ the 2-dimensional natural representation of $A_{1}$.  A basis for $\Sym{3}{V}$ is $\{ \smun{e}{1}^{3}, \smun{e}{1}^{2}\smun{e}{2}, \smun{e}{1}\smun{e}{2}^{2}, \smun{e}{2}^{3} \}$ and the highest weight vector is $\smun{e}{1}^{3}$.  The dimension of $\Sym{3}{V}$ is 4.  The corresponding $\smun{\Lie{g}}{0}$-module isomorphism is $\smun{e}{1}^{3} \mapsto \smun{X}{2}^{-}, \smun{e}{1}^{2}\smun{e}{2} \mapsto \smun{X}{12}^{-}, \smun{e}{1}\smun{e}{2}^{2} \mapsto \smun{X}{112}^{-}, \smun{e}{2}^{3} \mapsto \smun{X}{1112}^{-}$.

$\Lie{b}_{-2}$ has highest weight $\omega_0$ so is the trivial representation, by a dimension calculation.  It is spanned by the highest root in $G_{2}$, $\smun{X}{11122}^{-}$.  We can consider $\Lie{b}_{-2}$ to be spanned by $\smun{e}{1}^{3}\smun{e}{2}^{3}$, for the following reason.

The bracket on $\Lie{b} = \Lie{b}_{-1} \dsum \Lie{b}_{-2}$ is a cocycle central extension of the zero bracket on $\Lie{b}_{-1}$, given explicitly by $\Lbracket{\smun{e}{1}^{i}\smun{e}{2}^{j}}{\smun{e}{1}^{k}\smun{e}{2}^{l}}=\delta_{(i+k),3}\, \delta_{(j+l),3}\, \smun{e}{1}^{3}\smun{e}{2}^{3}$.  The braided-Lie cobracket is \[ \underline{\delta}( \smun{e}{1}^{3}\smun{e}{2}^{3} ) = \sum_{\substack{i,j,k,l=0 \\ i+k = 3 \\ j+l =3}}^{3} \gamma_{ijkl}\, \smun{e}{1}^{i}\smun{e}{2}^{j}\wedge \smun{e}{1}^{k}\smun{e}{2}^{l} \] for some non-zero constants $\gamma_{ijkl}$.

\end{description}

\begin{jegnote}{} This case has been covered as Example 4.6 in \cite{BraidedLie}. \end{jegnote}

\begin{description}
\item[($E_{7}D_{6}$)] Deletion $(E_{7}, 1, D_{6}, i \mapsto 8-i )$

$\Lie{b}_{-1}$ has highest weight $\omega_5$ so is the negative (half-) spin representation $S_{6}^{-}$ of $D_{6}$ (see for example \cite[Chapter 20]{FultonHarris}).  The dimension of $S_{6}^{-}$ is 32.  A basis for $S_{6}^{-}=V \dsum \Wedge{3}{(V)} \dsum \Wedge{5}{(V)}$ (as vector spaces; $V$ the vector space of dimension 6) is given by taking the natural bases for these pieces and the highest weight vector is $\smun{e}{1} \wedge \smun{e}{2} \wedge \smun{e}{3} \wedge \smun{e}{4} \wedge \smun{e}{5}$. 

$\Lie{b}_{-2}$ has highest weight $\omega_0$ so is the trivial representation.  In what follows, we see that $\Wedge{6}{(V)}$, spanned by $\smun{e}{1} \wedge \smun{e}{2} \wedge \smun{e}{3} \wedge \smun{e}{4} \wedge \smun{e}{5} \wedge \smun{e}{6}$, is the correct choice of basis for $\Lie{b}_{-2}$.

The bracket on $\Lie{b}=\Lie{b}_{-1} \dsum \Lie{b}_{-2}$ is given by the wedge product, \ie is non-zero on the subspaces $V \wedge \Wedge{5}{(V)}$ and $\Wedge{3}{(V)} \wedge \Wedge{3}{(V)}$ of $\Lie{b}_{-1} \wedge \Lie{b}_{-1}$.  The braided-Lie cobracket will be a non-zero map $\underline{\delta} : \Wedge{6}{(V)} \to \Wedge{6}{(V)}$, \ie is a non-zero scalar multiple of the identity.

\item[($BA$)] Deletion $(B_{l+1}, l+1, A_{l}, i \mapsto l-i+1 )$

$\Lie{b}_{-1}$ has highest weight $\omega_1$ so is the natural representation of $A_{l}$ on the vector space $V$ of dimension $l+1$.  A basis for $V$ is $\{ \smun{e}{i} \mid 1\leq i \leq l+1 \}$ and the highest weight vector is $\smun{e}{1}$.  The corresponding $\smun{\Lie{g}}{0}$-module isomorphism is $\smun{e}{i} \mapsto \smun{X}{i(i+1)\cdots (l+1)}^{-}$, for $1\leq i \leq l+1$.

$\Lie{b}_{-2}$ has highest weight $\omega_2$ so is the second exterior power $\Wedge{2}{(V)}$ with $V$ as before.  The dimension of $\Wedge{2}{(V)}$ is $\frac{1}{2}l(l+1)$.  A basis for $\Wedge{2}{(V)}$ is $\{ \smun{e}{i} \wedge \smun{e}{j} \mid 1\leq i < j \leq l+1 \}$ and the highest weight vector is $\smun{e}{1} \wedge \smun{e}{2}$.  We may deduce this from the following.

The bracket on $\Lie{b}=\Lie{b}_{-1} \dsum \Lie{b}_{-2}$ is non-zero: for example, there exists $X_{\alpha}^{-}$ such that $\Lbracket{ X_{\alpha}^{-}}{\smun{X}{l+1}^{-}} \in \Lie{g}^{\Lambda}$ where $\Lambda$ is the highest root in $\Lie{g}=B_{l+1}$.  Thus $\Lie{b}_{-2} \iso \Wedge{2}{\Lie{b}_{-1}} / \Ker{\Lbracket{\ }{\ }_{-1}}$ but $\Wedge{2}{\Lie{b}_{-1}}=\Wedge{2}{(V)}$ is irreducible.  Since $\Ker{\Lbracket{\ }{\ }_{-1}}\neq \Wedge{2}{\Lie{b}_{-1}}$, we see that $\Lie{b}_{-2} \iso \Wedge{2}{\Lie{b}_{-1}}=\Wedge{2}{(V)}$.  Further, $\Lbracket{\ }{\ }_{-1}=\wedge : V\tensor V \to \Wedge{2}{(V)}$.  The braided-Lie cobracket $\underline{\delta}: \Lie{b}_{-2} \to \Wedge{2}{\Lie{b}_{-1}}$ is an isomorphism.

The $\smun{\Lie{g}}{0}$-module isomorphism is given on $\Lie{b}_{-2}$ by  \[ \smun{e}{i} \wedge \smun{e}{j} \mapsto \smun{X}{i(i+1)\cdots (j-1)(j)(j)\cdots (l+1)(l+1)}^{-} \] for $1\leq i < j \leq l+1$.

\item[($E_{7}A_{6}$)] Deletion $(E_{7}, 2, A_{6}, \left( \begin{array}{llllll} \Ss 1 & \Ss 2 & \Ss 3 & \Ss 4 & \Ss 5 & \Ss 6 \\ \Ss 1 & \Ss 3 & \Ss 4 & \Ss 5 & \Ss 6 & \Ss 7 \end{array} \right) )$

$\Lie{b}_{-1}$ has highest weight $\omega_3$ so is the third exterior power $\Wedge{3}{(V)}$ with $V$ the 7-dimensional natural representation of $A_{6}$.  The dimension of $\Wedge{3}{(V)}$ is 35.  A basis for $V$ is $\{ \smun{e}{i} \wedge \smun{e}{j} \smun{e}{k} \mid 1\leq i < j < k \leq 7 \}$ and the highest weight vector is $\smun{e}{1} \wedge \smun{e}{2} \wedge \smun{e}{3}$. 

$\Lie{b}_{-2}$ has highest weight $\omega_6$ by considering the module decomposition $\Wedge{2}{\Lie{b}_{-1}}=V([0,1,0,1,0,0]) \dsum \Wedge{6}{(V)}$ (we use a formula in \cite[Chapter 15]{FultonHarris}) and a dimension calculation.  We use the usual natural basis for $\Wedge{6}{(V)}$ rather than a basis in terms of the dual of $V$, even though $\Wedge{6}{(V)} \iso \dual{V}$.  The dimension of $\Wedge{6}{(V)}$ is 7.

The bracket on $\Lie{b}=\Lie{b}_{-1} \dsum \Lie{b}_{-2}$ is given by the wedge product map $\wedge : \Lie{b}_{-1} \tensor \Lie{b}_{-1} = \Wedge{3}{(V)} \tensor \Wedge{3}{(V)} \to \Lie{b}_{-2} = \Wedge{6}{(V)}$.  The kernel of $\wedge$ is $V([0,1,0,1,0,0])$.  The braided-Lie cobracket $\underline{\delta}$ is an isomorphism.
\item[($F_{4}B_{3}$)] Deletion $(F_{4}, 4, B_{3}, \id )$

$\Lie{b}_{-1}$ has highest weight $\omega_3$ so is the 8-dimensional spinor representation $S_{3}$ of $B_{3}=\Lie{so}_{7}$.  A basis for $S_{3}=\bigdsum_{i=0}^{3} \Wedge{i}{(V)}$ (as vector spaces; $V$ the vector space of dimension 3) is given by taking the natural basis for each piece and the highest weight vector is $\smun{e}{1} \wedge \smun{e}{2} \wedge \smun{e}{3}$.

$\Lie{b}_{-2}$ has highest weight $\omega_1$, by considering the module decomposition $\Wedge{2}{\Lie{b}_{-1}} = \Wedge{2}{(S_3)} \iso W \dsum \Wedge{2}{(W)}$ for $W$ the 7-dimensional natural representation of $B_3$ and a dimension calculation.  We obtain this decomposition by examining the above description of $S_3$.  So, $\Lie{b}_{-2}$ is isomorphic to the natural representation, $W$.

The bracket on $\Lie{b}=\Lie{b}_{-1} \dsum \Lie{b}_{-2}$ does not seem to have an interpretation as a natural map on $\Wedge{2}{(S_3)}$.

\item[($E_{8}D_{7}$)] Deletion $(E_{8}, 1, D_{7}, i\mapsto 9-i )$

$\Lie{b}_{-1}$ has highest weight $\omega_6$ so is the positive (half-)spin representation $S_{7}^{+}$ of $D_{7}=\Lie{so}_{14}$.  The dimension of $S_{7}^{+}$ is 64.  As a vector space, we have $S_{7}^{+}=\Wedge{\text{even}}{(V)}=\bigdsum_{i=0,2,4,6} \Wedge{i}{(V)}$ with $V$ the vector space of dimension 7 so a basis is given by taking the natural basis for each piece. The highest weight vector is $\smun{e}{1} \wedge \smun{e}{2} \wedge \smun{e}{3} \wedge \smun{e}{4} \wedge \smun{e}{5} \wedge \smun{e}{6}$.

$\Lie{b}_{-2}$ has highest weight $\omega_1$, so is the 14-dimensional natural representation $W$ of $D_7$.  A basis for $W$ is $\{ \smun{e}{i} \mid 1 \leq i \leq 14 \}$ and the highest weight vector is $\smun{e}{1}$.  We obtain this from the decomposition $\Wedge{2}{\Lie{b}_{-1}}=\Wedge{2}{(S_{7}^{+})} \iso \Wedge{5}{(W)} \dsum W$ and a dimension calculation.
\end{description}

\subsection{$m_{d}=3$}\label{md3}

\begin{description}
\item[($E_{8}A_{7}$)] Deletion $(E_{8}, 2, A_{7}, \left( \begin{array}{lllllll} \Ss 1 & \Ss 2 & \Ss 3 & \Ss 4 & \Ss 5 & \Ss 6 & \Ss 7 \\ \Ss 1 & \Ss 3 & \Ss 4 & \Ss 5 & \Ss 6 & \Ss 7 & \Ss 8 \end{array} \right) )$

$\Lie{b}_{-1}$ has highest weight $\omega_3$ so is the third exterior power $\Wedge{3}{(V)}$ with $V$ the 8-dimensional natural representation of $A_{7}$.  The dimension of $\Wedge{3}{(V)}$ is 56.  We take the natural basis for $\Wedge{3}{(V)}$ and the highest weight vector is $\smun{e}{1} \wedge \smun{e}{2} \wedge \smun{e}{3}$.

$\Lie{b}_{-2}$ has highest weight $\omega_6$ so is the sixth exterior power $\Wedge{6}{(V)}$, with $V$ as before.  The dimension of $\Wedge{6}{(V)}$ is 28.  We take the natural basis and the highest weight vector is $\smun{e}{1} \wedge \cdots \wedge \smun{e}{6}$.  We obtain this from the decomposition $\Wedge{2}{\Lie{b}_{-1}} = \Wedge{2}{(\Wedge{3}{(V)})}\iso V([0,1,0,1,0,0,0]) \dsum \Wedge{6}{(V)}$ and calculating dimensions.
 
$\Lie{b}_{-3}$ has highest weight $\omega_1$, so is the 8-dimensional natural representation $V$.  The highest weight vector is $\smun{e}{1}$.  We see this since the tensor product of $\Lie{b}_{-1} \tensor \Lie{b}_{-2}$ decomposes as \[ \Wedge{3}{(V)} \tensor \Wedge{6}{(V)} \iso V([0,0,1,0,0,1,0]) \dsum V([0,1,0,0,0,0,1]) \dsum V \] so, by the same arguments about the kernel of the bracket map, we can use a dimension calculation as before.

\item[($G_{2}A_{1}$)(b)] Deletion $(G_{2}, 1, A_{1}, \left( \begin{array}{l} \Ss 1 \\ \Ss 2 \end{array} \right) )$

$\Lie{b}_{-1}$ has highest weight $\omega_1$ so is the 2-dimensional natural representation of $A_{1}=\Lie{sl}_{2}$.  A basis for $V$ is $\{ \smun{e}{1}, \smun{e}{2} \}$ and the highest weight vector is $\smun{e}{1}$.  We have $\smun{e}{1} \mapsto \smun{X}{1}^{-},\ \smun{e}{2} \mapsto \smun{X}{12}^{-}$.

$\Lie{b}_{-2}$ has highest weight $\omega_0$ so is the trivial representation, spanned by $\varsigma$, say.  This is since $\Wedge{2}{(V)}\iso \complex$ and $\Ker{\Lbracket{\ }{\ }_{-1}}=0$ (the bracket is non-zero on $\Lie{b}_{-1}$).  We have $\varsigma \mapsto \smun{X}{112}^{-}$.

$\Lie{b}_{-3}$ has highest weight $\omega_{1}$ so is another copy of the natural representation $V$, with basis $\{ \smun{f}{1}, \smun{f}{2} \}$.  The highest weight vector is $\smun{f}{1}$.  This is obtained from a direct examination of the root system of $G_2$, giving $\smun{f}{1} \mapsto \smun{X}{1112}^{-}$ and $\smun{f}{2} \mapsto \smun{X}{11122}^{-}$.  

The bracket in these bases is \begin{align*} \Lbracket{\smun{e}{1}}{\smun{e}{2}}& =c_{1}\,\varsigma, \\ \Lbracket{\smun{e}{1}}{\varsigma} & = c_{2}\, \smun{f}{1}, \\ \Lbracket{\smun{e}{2}}{\varsigma} & = c_{3}\, \smun{f}{2} \end{align*} for some constants $c_i$.

\end{description}

\bibliographystyle{halpha}
\bibliography{references}

\vspace{1em}

\noindent\begin{tabular}{ll}
Address:  & 	School of Mathematical Sciences, \\
	  & 	Queen Mary, University of London, \\
	  &	E1 4NS, \\
	  &	United Kingdom. \\
	  &	\\
E-mail:   &	J.Grabowski@qmul.ac.uk \\
Web site: &     http://www.maths.qmul.ac.uk/\webtilde jeg/ \\
	  &	\\
MSC:	  &	17Bxx (Primary), 22Exx
\end{tabular}

\end{document}